# Discontinuous Petrov-Galerkin boundary elements [*]

Norbert Heuer[†]    Michael Karkulik[†]

**Abstract**

Generalizing the framework of an ultra-weak formulation for a hypersingular integral equation on closed polygons in [N. Heuer, F. Pinochet, arXiv 1309.1697 (to appear in *SIAM J. Numer. Anal.*)], we study the case of a hypersingular integral equation on open and closed polyhedral surfaces. We develop a general ultra-weak setting in fractional-order Sobolev spaces and prove its well-posedness and equivalence with the traditional formulation. Based on the ultra-weak formulation, we establish a discontinuous Petrov-Galerkin method with optimal test functions and prove its quasi-optimal convergence in related Sobolev norms. For closed surfaces, this general result implies quasi-optimal convergence in the $L^2$-norm. Some numerical experiments confirm expected convergence rates.

## 1 Introduction

In recent years, the discontinuous Petrov-Galerkin (DPG) method with optimal test functions has drawn some attention. In most cases it is based on an ultra-weak formulation (cf. Després and Cessenat [10, 19]) that one obtains by integrating by parts element-wise the underlying first order system and by replacing boundary terms with new unknowns (cf. Botasso, Micheletti and Sacco [4]). Choosing appropriate test functions, discrete stability can be deduced from the stability of the continuous formulation. In this combination, the method has been proposed and analyzed by Demkowicz and Gopalakrishnan, cf. [15, 16], with particular emphasis on its possible robustness for singularly perturbed problems [5, 7, 11, 17, 18].

In this paper, we develop a DPG method with optimal test functions to numerically solve a hypersingular integral equation on (open or closed) polyhedral surfaces. Our model problem governs the Laplacian and there is no reason to be concerned about stability (the standard Galerkin method is stable). However, as for partial differential equations, our hypothesis is that this strategy can lead to stable approximations for singularly perturbed integral equations stemming, e.g., from acoustic scattering or almost incompressible linear elasticity. We refer to [9, 12, 13, 20, 37] for some discussions of stability issues in the boundary element Galerkin method. We follow the abstract framework from [31] where a DPG method for a hypersingular integral equation has been studied in the two-dimensional case of a closed polygonal curve. Here, we deal with the three-dimensional polyhedral case and consider, in particular, open surfaces. Open curves were excluded in [31].

---

[*]Supported by CONICYT through FONDECYT projects 1110324, 3140614 and Anillo ACT1118 (ANANUM).
[†]Facultad de Matemáticas, Pontificia Universidad Católica de Chile, Avenida Vicuña Mackenna 4860, Macul, Santiago, Chile, email: {nheuer,mkarkulik}@mat.puc.cl

The situation of a closed surface is most convenient since, in that case, the whole formulation is set in simplest Sobolev spaces. The solution is approximated in $L^2$ and test functions are taken from $L^2$ and piecewise $H^1$-spaces. In this way (apart from a skeleton variable in $H^{-1/2}$, which is not being controlled) we completely avoid fractional-order Sobolev spaces. For the first time, a hypersingular integral equation in three dimensions is weakly formulated and approximated in $L^2$ and $H^1$-spaces without relying on complicated dualities usually induced by elliptic operators of order one. The situation on an open surface is slightly more involved. This is due to the fact that, in this case, solutions have strong edge singularities so that corresponding $H^1$-estimates (needed for the analysis of the dual problem) must be avoided. Nevertheless, also in this extreme case we are able to present a well-posed formulation and prove quasi-optimal convergence of the DPG method in standard Sobolev spaces. These spaces are of fractional order, but close to $L^2$ and $H^1$ so that no variational crimes arise.

Apart from the new mathematical framework that provides optimal results in standard Sobolev norms, there are practical advantages of the DPG method which apply also in the case of hypersingular operators. We recall the list from [31].

- The matrices of linear systems used for the approximation of optimal test functions and for error calculation are sparse.

- System matrices are symmetric and positive definite.

- Error control is inherent since errors in the energy norm can be calculated through the implementation of the trial-to-test operator.

- Since norms are localizable, the energy norm of the error gives local information which can be used to steer adaptive refinements.

- Error estimates and stability hold for any combination of meshes and polynomial degrees so that $hp$ methods do not require a new analysis.

- Since approximation spaces can be discontinuous, one has full flexibility for $h$ and $p$ adaptivity.

For a detailed discussion of these facts we refer to the previously mentioned references on the DPG method with optimal test functions.

There are also limitations of the DPG method. Most importantly, despite of using localizable norms, stiffness matrices are still densely populated when discretizing non-local operators. The global coupling of unknowns enters through the trial-to-test operator since it carries the information of the original operator. Secondly, optimal test functions must be approximated (except for particular cases, see [31]). The influence of this approximation on stability and convergence estimates has been analyzed for partial differential equations, see [24, 30], but is an open problem when dealing with integral operators. Finally, strict implementation of our method on open surfaces requires the use of inner products from (slightly) fractional-order Sobolev spaces (namely, of orders $-\epsilon$ and $1-\epsilon$ for fixed $\epsilon > 0$). This is not acceptable in practice and one simply selects $\epsilon = 0$ (corresponding to the parameter $s = 1/2$ in our estimates). In this limit case we still have quasi-optimal convergence, but the error is controlled in the energy norm defined by the problem rather than in $L^2$ (see Corollary 18 at the end of Section 4 for this result). Some more comments are made in Remark 19, after the corollary.

For convenience of the reader, we briefly recall the abstract framework of the DPG method with optimal test functions. For details and proofs we refer to [16, 39]. For Hilbert spaces $U$, $V$, bilinear form $b: U \times V \to \mathbb{R}$ and $L \in V'$, let

$$u \in U: \quad b(u, v) = L(v) \quad \forall v \in V$$

be a well-posed variational formulation. For an approximation space $U_{hp} \subset U$, the Petrov-Galerkin method with optimal test functions consists in calculating $u_{hp} \in U_{hp}$ such that

$$b(u_{hp}, v) = L(v) \quad \forall v \in \Theta(U_{hp}).$$

Here, $\Theta: U \to V$ is the trial-to-test operator defined by

$$\langle \Theta \varphi, v \rangle_V = b(\varphi, v) \quad \forall v \in V$$

with inner product $\langle \cdot, \cdot \rangle_V$ in $V$. Standard arguments from functional analysis show that $u_{hp}$ is the best approximation of $u$ in the energy norm defined by

$$\|\varphi\|_U := \sup_{v \in V \setminus \{0\}} \frac{b(\varphi, v)}{\|v\|_V}.$$

Now, practicality and efficiency of the method hinge on several ingredients.

- The trial-to-test operator $\Theta$ is localized by using discontinuous finite elements (the Petrov-Galerkin method is then called discontinuous Petrov-Galerkin method, DPG). Of course, this goes in hand with an appropriate bilinear form $b(\cdot, \cdot)$.

- One wants to select an appropriate norm in $U$ to control the error. For a given $V$-norm, the energy norm $\|\cdot\|_U$ is not necessarily convenient. Using a different norm in $U$, duality via $b(\cdot, \cdot)$ induces a norm in $V$ that is different from the inner-product norm.

- Efficient implementation of $\Theta$ not only requires discontinuous finite elements but also a localizable inner product in $V$. This conflicts with the selection of a norm in $U$.

- Even a localized version of $\Theta$ cannot be exactly implemented since local spaces are still infinite-dimensional. Therefore, in practice the trial-to-test operator is approximated by selecting finite-dimensional subspaces $V^r \subset V$ and by defining an approximating operator $\Theta^r : U \to V^r$ by $\langle \Theta^r \varphi, v \rangle_V = b(\varphi, v)$ for all $v \in V^r$. (This issue is not further analyzed in this paper. We refer to the previously mentioned references [24, 30] that deal with partial differential equations.)

Principal focus of DPG analysis is to deal with conflicting selections of norms in $U$ and $V$. The aim is to choose a norm in $U$ that is appropriate for the problem under consideration, and to find a corresponding $V$-norm that (i) allows for an efficient implementation of the DPG method and (ii) that guarantees a robust error estimate with appropriate convergence order.

In this paper, we propose and analyze a DPG method with optimal test functions for a hypersingular integral equation on polyhedral surfaces. First principal step is to develop a well-posed ultra-weak variational formulation. This is done in Section 4.1. In Section 4.2 we analyze the equivalence of different norms in the ansatz space $U$ and the test space $V$. Solvability (Theorem 14) and stability (Theorem 15) of the ultra-weak formulation is proved in Section 4.3. In Section 4.4 we present the DPG method and our main result (Theorem 17) which is a general Céa estimate. All these results depend on a parameter $s$ that serves to select specific norms in $U$ and $V$, the most convenient case being $s = 1/2$. However, for an open surface there are some complications with $s = 1/2$. This situation is being considered in Corollaries 16 and 18.

The remaining parts of this paper are as follows. In the next section we introduce the model problem and define partitions of the surface. For the convenience of the reader, in Section 2.1 we resume the simplest case (parameter $s = 1/2$) for a closed surface and present the corresponding main result (Theorem 1) which is a particular case of Theorem 17 from Section 4. In Section 3 we collect all the technical results from Sobolev spaces that are not directly related to the DPG analysis. This includes properties of integral operators and regularity of solutions. The principal part is Section 4, as discussed above. Some numerical experiments are presented in the last section.

Throughout the paper, $a \lesssim b$ means that $a \leq cb$ with a generic constant $c > 0$ that is independent of involved parameters and functions. Similarly, we use the notation $a \gtrsim b$ and $a \simeq b$.

## 2  Model problem and the closed-surface case

Let $\Gamma$ be the boundary of a simply-connected polyhedral domain $\Omega$, or a connected union of some of its faces. Our model problem is the hypersingular integral equation on $\Gamma$,

$$\mathcal{W}\phi(x) := -\frac{1}{4\pi} \partial_{\mathbf{n}(x)} \int_\Gamma \phi(y) \partial_{\mathbf{n}(y)} \frac{1}{|x-y|} \, ds_y = f(x) \quad (x \in \Gamma). \tag{1}$$

Here, $\mathbf{n}$ is the exterior unit normal vector on $\partial \Omega$, and $f$ is a given function. In the case that $\Gamma = \partial \Omega$, (1) represents the Neumann problem for the Laplacian in $\mathbb{R}^3 \setminus \overline{\Omega}$, when selecting $f = (1/2 - K')v$ with Neumann datum $v$ and $K'$ being the adjoint of the so-called double-layer operator, cf. [32].

A weak form of the hypersingular integral equation is

$$\phi \in \widetilde{H}^{1/2}(\Gamma): \quad \langle \mathcal{W}\phi, \psi \rangle_\Gamma + m_\Gamma(\phi, \psi) = \langle f, \psi \rangle_\Gamma \quad \forall \psi \in \widetilde{H}^{1/2}(\Gamma) \tag{2}$$

where

$$m_\Gamma(\phi, \psi) := \begin{cases} 0 & \text{if } \Gamma \text{ is an open surface,} \\ \langle \phi, 1 \rangle_\Gamma \langle \psi, 1 \rangle_\Gamma & \text{if } \Gamma \text{ is closed.} \end{cases} \tag{3}$$

The space $\widetilde{H}^{1/2}(\Gamma)$ is the trace space of $H^1(\Omega)$ when $\Gamma = \partial\Omega$, and consists of trace functions that vanish on $\partial\Omega \setminus \Gamma$ when $\Gamma \neq \partial\Omega$. Furthermore, $\langle \cdot, \cdot \rangle_\Gamma$ denotes the $L^2(\Gamma)$-bilinear form and its extension by duality. The rank-one term $m_\Gamma(\cdot, \cdot)$ eliminates the kernel of $\mathcal{W}$ in the case of a closed surface.

In the following, we also need the single layer operator $\mathcal{V}$ defined by

$$\mathcal{V}\phi(x) := \frac{1}{4\pi} \int_\Gamma \frac{\phi(y)}{|x-y|} \, ds_y \quad (x \in \Gamma),$$

and surface differential operators curl (scalar) and **curl** (vectorial), defined as follows. Let $\nabla$ denote the surface gradient on $\partial\Omega$. Following [8], the surface vector curl operator $\mathbf{curl}\, v = -\mathbf{n} \times \nabla v$ can be defined on $H^{1/2}(\partial\Omega)$. It maps to a space of tangential vector fields. By curl, we denote the adjoint operator of **curl** with respect to the $L^2$-bilinear form. On a smooth surface it is given by $\text{curl}\, \boldsymbol{\sigma} = -\nabla \cdot (\mathbf{n} \times \boldsymbol{\sigma})$, where $\nabla \cdot$ denotes the surface divergence.

An ultra-weak formulation of (1) and resulting DPG method is based on piecewise integration by parts, i.e., it hinges on a partition of $\Gamma$. Let $\mathcal{T}$ denote such a partition (also called mesh) that is compatible with the geometry, i.e., $\mathcal{T}$ is a finite set, the elements $T \in \mathcal{T}$ are mutually disjoint, open sets with $\bigcup_{T \in \mathcal{T}} \overline{T} = \overline{\Gamma}$, and every element is a subset of a face of $\Gamma$. Furthermore, for simplicity, we assume that every $T \in \mathcal{T}$ is an affine image of one of a finite set of reference elements $\widehat{T} \subset \mathbb{R}^2$. As $\mathcal{T}$ is finite, all its elements $T$ are shape-regular in the sense that (with $\omega = T$)

$$\gamma_\omega := \frac{\text{diam}(\omega)}{\sup\{\text{diam}(B);\ B \subset \omega \text{ is a (hyper) circle}\}} < \infty.$$

Let us denote by $\gamma_\mathcal{T} := \max\{\gamma_T;\ T \in \mathcal{T}\}$ the shape-regularity parameter of $\mathcal{T}$. Throughout we assume that $\gamma_\mathcal{T}$ is bounded. We also introduce a mesh parameter by $h := \min\{h_T;\ T \in \mathcal{T}\}$ with $h_\omega := (\int_\omega 1\, ds)^{1/2}$ for any $\omega \subset \Gamma$.

In the remainder of the paper we establish a general ultra-weak formulation of the hypersingular integral equation, prove its well-posedness, and analyze a DPG method that is based on the new formulation. In the case of an open surface, details are quite technical and require some tedious definitions of particular norms. In the following section, for ease of reading, we therefore present the case of a closed surface which can be much simplified. In particular, we propose an ultra-weak formulation, present the DPG method with optimal test functions, and announce one of our main results that establishes quasi-optimal convergence of the error in $L^2$-norms. The analysis of this situation is simply a particular case in Section 4.

## 2.1 Ultra-weak formulation and DPG method for a closed surface

We need to define some Sobolev spaces which will be generalized later in this paper. We use the standard $L^2$ and $H^1$ notation, and use bold symbols to indicate (tangential) vector-valued functions and their spaces. Let us mention at this point that there are three components of the ultra-weak solution, $\phi$, $\boldsymbol{\sigma} = \mathcal{V}\mathbf{curl}\,\phi$, and tangential components of $\boldsymbol{\sigma}$ on the skeleton of $\mathcal{T}$. Both $\phi$ and $\boldsymbol{\sigma}$ will be measured in $L^2$, and for the tangential components of $\boldsymbol{\sigma}$ we need the interface space

$$H^{-1/2}(\mathcal{S}) := \Big\{ \hat{\sigma} \in \prod_{T \in \mathcal{T}} H^{-1/2}(\partial T); \, \exists \boldsymbol{\sigma} \in \mathbf{H}(\mathrm{curl},\Gamma) \text{ such that } \hat{\sigma}|_{\partial T} = \boldsymbol{\sigma} \cdot \boldsymbol{t}|_{\partial T} \, \forall T \in \mathcal{T} \Big\}.$$

Here, the symbol $\mathcal{S}$ denotes the skeleton of the partition $\mathcal{T}$ (the collection of boundaries $\partial T$ with $T \in \mathcal{T}$) and $\mathbf{H}(\mathrm{curl},\Gamma)$ consists of $L^2(\Gamma)$-functions $\boldsymbol{\sigma}$ with $\mathrm{curl}\,\boldsymbol{\sigma} \in L^2(\Gamma)$, equipped with the standard norm. On the boundary $\partial T$ of an element $T \in \mathcal{T}$, $\boldsymbol{t}$ is the unit tangential vector along $\partial T$ and $\boldsymbol{\sigma} \cdot \boldsymbol{t}|_{\partial T}$ is the canonical tangential trace of a function from $\mathbf{H}(\mathrm{curl},T)$. The space $H^{-1/2}(\mathcal{S})$ is normed by

$$\|\hat{\sigma}\|_{H^{-1/2}(\mathcal{S})} := \inf_{\boldsymbol{\sigma} \cdot \boldsymbol{t}|_{\mathcal{S}} = \hat{\sigma}} \|\boldsymbol{\sigma}\|_{\mathbf{H}(\mathrm{curl},\Gamma)}.$$

Test functions will be taken from $\mathbf{L}^2(\Gamma)$ and the piecewise $H^1$-space

$$H^1(\mathcal{T}) := \big\{ v \in L^2(\Gamma); \, v|_T \in H^1(T) \, \forall T \in \mathcal{T} \big\}$$

with canonical product norm.

**Ultra-weak formulation.** For given $f \in L^2(\Gamma)$ with $\langle f, 1 \rangle_\Gamma = 0$, the ultra-weak formulation of (2) reads as follows: *Find $(\boldsymbol{\sigma}, \phi, \hat{\sigma}) \in U^{1/2} := \mathbf{L}^2(\Gamma) \times L^2(\Gamma) \times H^{-1/2}(\mathcal{S})$ such that*

$$\langle \boldsymbol{\sigma}, \boldsymbol{\tau} \rangle_\Gamma - \langle \phi, \mathrm{curl}\,\mathcal{V}\boldsymbol{\tau} \rangle_\Gamma = 0,$$
$$\langle \boldsymbol{\sigma}, \mathbf{curl}_{\mathcal{T}} v \rangle_\Gamma + \langle \hat{\sigma}, [v] \rangle_{\mathcal{S}} + m_\Gamma(\phi, v) = \langle f, v \rangle_\Gamma$$

*for all $(\boldsymbol{\tau}, v) \in V^{1/2} := \mathbf{L}^2(\Gamma) \times H^1(\mathcal{T})$.* Here, $\mathbf{curl}_{\mathcal{T}}$ is the piecewise $\mathbf{curl}$-operator, and $[v]$ indicates the jump of $v$ across $\mathcal{S}$, appearing through the definition

$$\langle \hat{\sigma}, [v] \rangle_{\mathcal{S}} := \sum_{T \in \mathcal{T}} \langle \hat{\sigma}|_T, v|_T \rangle_{\partial T}.$$

Note that $\hat{\sigma}|_{\partial T} = \boldsymbol{\sigma} \cdot \boldsymbol{t}|_{\partial T}$ and that the tangential vector $\boldsymbol{t}$ on an edge of the mesh (which does not lie on $\partial \Gamma$) has different directions on the elements that share the edge. The term $m_\Gamma$ ensures that $\int_\Gamma \phi = 0$. We abbreviate the ultra-weak formulation as:

$$(\boldsymbol{\sigma}, \phi, \hat{\sigma}) \in U^{1/2}: \quad b(\boldsymbol{\sigma}, \phi, \hat{\sigma}; \boldsymbol{\tau}, v) = \langle f, v \rangle_\Gamma \quad \forall (\boldsymbol{\tau}, v) \in V^{1/2}, \tag{4}$$

with bilinear form $b: U^{1/2} \times V^{1/2} \to \mathbb{R}$ defined by

$$b(\boldsymbol{\sigma}, \phi, \hat{\sigma}; \boldsymbol{\tau}, v) := \langle \boldsymbol{\sigma}, \boldsymbol{\tau} + \mathbf{curl}_{\mathcal{T}} v \rangle_\Gamma - \langle \phi, \mathrm{curl}\,\mathcal{V}\boldsymbol{\tau} \rangle_\Gamma + \langle \hat{\sigma}, [v] \rangle_{\mathcal{S}} + m_\Gamma(\phi, v).$$

The spaces $U^{1/2}$ and $V^{1/2}$ are Hilbert spaces with canonical norms

$$\|(\boldsymbol{\sigma},\phi,\hat{\sigma})\|_{U^{1/2}} := \left(\|\boldsymbol{\sigma}\|_{\mathbf{L}^2(\Gamma)}^2 + \|\phi\|_{L^2(\Gamma)}^2 + \|\hat{\sigma}\|_{H^{-1/2}(\mathcal{S})}^2\right)^{1/2}$$

and

$$\|(\boldsymbol{\tau},v)\|_{V^{1/2}} := \left(\|\boldsymbol{\tau}\|_{\mathbf{L}^2(\Gamma)}^2 + \|v\|_{H^1(\mathcal{T})}^2\right)^{1/2}.$$

**DPG method.** Consider a discrete space $U_{hp} \subset U^{1/2}$ (usually piecewise polynomial functions with respect to the mesh $\mathcal{T}$) and let

$$\{\boldsymbol{u}_i;\ i=1,\ldots,\dim(U_{hp})\}$$

be a basis of $U_{hp}$. Then we define the discrete test space $V_{hp} \subset V^{1/2}$ with basis given by

$$\{\mathbf{v}_i := \Theta(\boldsymbol{u}_i);\ i=1,\ldots,\dim(U_{hp})\}.$$

Here, $\Theta : U^{1/2} \to V^{1/2}$ is the trial-to-test operator defined by

$$\langle\Theta(\boldsymbol{u}_i),\mathbf{v}\rangle_{V^{1/2}} = b(\boldsymbol{u}_i,\mathbf{v}) \quad \forall \mathbf{v} \in V^{1/2}$$

with $V^{1/2}$-inner product $\langle\cdot,\cdot\rangle_{V^{1/2}}$ that induces the norm $\|\cdot\|_{V^{1/2}}$. The DPG method with optimal test functions consists in finding $(\boldsymbol{\sigma}_{hp},\phi_{hp},\hat{\sigma}_{hp}) \in U_{hp}$ such that

$$b(\boldsymbol{\sigma}_{hp},\phi_{hp},\hat{\sigma}_{hp};\boldsymbol{\tau},v) = \langle f,v\rangle_\Gamma \quad \forall(\boldsymbol{\tau},v) \in V_{hp}. \tag{5}$$

There holds the following result.

**Theorem 1.** *For $f \in L^2(\Gamma)$ with $\langle f,1\rangle_\Gamma = 0$, the weak formulation (2) and ultra-weak formulation (4) are uniquely solvable and equivalent. The solution of (4) is given by $\boldsymbol{\sigma} = \mathcal{V}\mathbf{curl}\,\phi$, $\hat{\sigma} = \prod_{T\in\mathcal{T}} \boldsymbol{\sigma}\cdot\boldsymbol{t}|_{\partial T}$ and with $\phi$ being the solution of (2). Furthermore, there is a unique solution $(\boldsymbol{\sigma}_{hp},\phi_{hp},\hat{\sigma}_{hp}) \in U_{hp}$ of the DPG scheme (5). It satisfies the quasi-optimal error estimate*

$$\|\boldsymbol{\sigma}-\boldsymbol{\sigma}_{hp}\|_{\mathbf{L}^2(\Gamma)} + \|\phi-\phi_{hp}\|_{L^2(\Gamma)}$$
$$\lesssim \inf_{(\boldsymbol{\sigma}_d,\phi_d,\hat{\sigma}_d)\in U_{hp}} \left(\|\boldsymbol{\sigma}-\boldsymbol{\sigma}_d\|_{\mathbf{L}^2(\Gamma)} + \|\phi-\phi_d\|_{L^2(\Gamma)} + \|\hat{\sigma}-\hat{\sigma}_d\|_{H^{-1/2}(\mathcal{S})}\right).$$

*The hidden constant in the estimate is independent of $\mathcal{T}$ as long as the shape-regularity parameter $\gamma_\mathcal{T}$ is bounded.*

In what remains, a proof of this result is developed in a more general setting, comprising open and closed surfaces. The particular case of the theorem, for a closed surface, is obtained by combining the results of Theorems 14, 15, and 17 from Section 4 for the selection $s = 1/2$. Note that the regularity parameter $C_\Gamma(1/2)$ appearing in Theorems 15 and 17 is bounded in this case.

# 3 Basic results from Sobolev spaces

In this section we collect several technical results that will be needed for the analysis of our ultra-weak formulation and the DPG method. In the first part we give definitions of Sobolev spaces and norms, and recall properties of boundary integral operators. Furthermore, we introduce a regularity parameter $C_\Gamma(\cdot)$ for the solutions of hypersingular integral equations. This parameter enters the stability estimate of the general ultra-weak formulation and also appears in the quasi-optimal error estimate of the DPG method. In the second part, Section 3.2, we present equivalences and bounds for several fractional-order Sobolev norms and recall some properties of surface differential operators, including an integration-by-parts formula. We also introduce some further Sobolev spaces that are related with the surface differential operators.

## 3.1 Sobolev spaces, properties of integral operators, and regularity

From now on, $\omega \subset \partial\Omega$ always denotes a subset of the boundary $\partial\Omega$ of a bounded polyhedral domain $\Omega \subset \mathbb{R}^3$ such that $\omega$ itself has a Lipschitz boundary $\partial\omega$. The spaces $L^2(\omega)$ and $H^1(\omega)$ with norms $\|\cdot\|_{L^2(\omega)}$, $\|\cdot\|_{H^1(\omega)}$ and semi-norm $|\cdot|_{H^1(\omega)}$ denote usual Sobolev spaces. Vector valued versions of these spaces and their elements will be denoted by bold symbols, e.g., $\boldsymbol{\sigma} \in \mathbf{H}^1(\omega)$. More precisely, all the vector-valued functions and their spaces are tangential to the surface under consideration. The space $\widetilde H^1(\omega)$ (traditionally denoted by $H^1_0(\omega)$) is the space of $H^1(\omega)$ functions with vanishing trace on the boundary $\partial\omega$ equipped with the norm $|\cdot|_{H^1(\omega)}$. For $s \in (0,1)$, (semi-) norms are defined by

$$|v|^2_{H^s(\omega)} := \int_\omega \int_\omega \frac{|v(x)-v(y)|^2}{|x-y|^{2+2s}}\,ds_x\,ds_y, \text{ and}$$

$$\|v\|^2_{\widetilde H^s(\omega)} := \int_\omega \int_\omega \frac{|v(x)-v(y)|^2}{|x-y|^{2+2s}}\,ds_x\,ds_y + \int_\omega v(x)^2 \mathrm{dist}_\Gamma(x,\partial\omega)^{-2s}\,ds_x.$$

For $s \in (0,1)$, the space $H^s(\omega)$ is defined as the space of functions $v \in L^2(\omega)$ with $|v|_{H^s(\omega)} < \infty$ and norm $\|\cdot\|^2_{H^s(\omega)} := \|\cdot\|^2_{L^2(\omega)} + |\cdot|^2_{H^s(\omega)}$, while the space $\widetilde H^s(\omega)$ is defined as the space of functions $v \in L^2(\omega)$ with $\|v\|_{\widetilde H^s(\omega)} < \infty$. In some proofs, we will also employ Sobolev norms defined by interpolation, e.g., $H^s(\omega) = [L^2(\omega), H^1(\omega)]_s$ and $\widetilde H^s(\omega) = [L^2(\omega), \widetilde H^1(\omega)]_s$, $s \in (0,1)$, defined by the real K-method [3]. It is well known that their norms are equivalent with the corresponding norms of Sobolev-Slobodeckij type defined previously. Though equivalence numbers depend on $s$ and the domain $\omega$, cf. [28, 29]. Abusing notation, we set $\widetilde H^s(\Gamma) := H^s(\Gamma)$ with norm $\|\cdot\|_{H^s(\Gamma)}$ for a closed surface $\Gamma$. The space $H^s_h(\omega)$ for $s \in [0,1]$ is defined as $H^s(\omega)$ with norm $\|\cdot\|^2_{H^s_h(\omega)} := h_\omega^{-2s}\|\cdot\|^2_{L^2(\omega)} + |\cdot|^2_{H^s(\omega)}$. Recall the notation $h_\omega = (\int_\omega 1\,ds)^{1/2}$; it is a measure for the diameter of $\omega$. If a space carries the lower index 0, e.g., $H^s_0(\Gamma)$ and $L^2_0(\Gamma)$, it consists of functions with vanishing integral mean in case of a closed surface $\Gamma$. If $\Gamma$ is open the index 0 has no meaning.

Spaces with negative order $s < 0$, as well as their respective norms, are defined dual to spaces with positive order $-s > 0$. Duality is understood with respect to the extended $L^2$-inner product, which is denoted $\langle \cdot, \cdot \rangle$ with appropriate index to indicate the geometric object under consideration. For example, the dual space of $H^s(\omega)$ is $\widetilde{H}^{-s}(\omega)$, and the norm is given by

$$\|\phi\|_{\widetilde{H}^{-s}(\omega)} := \sup_{v \in H^s(\omega)} \frac{\langle \phi, v \rangle_\Gamma}{\|v\|_{H^s(\omega)}}.$$

Note that the duals to $H^s$ with $s > 0$ are denoted $\widetilde{H}^{-s}$, while the duals to $\widetilde{H}^s$ with $s > 0$ are denoted $H^{-s}$.

Product spaces are defined for $s \in [0,1]$ via

$$H^s(\mathcal{T}) := \{v \in L^2(\Gamma);\ v|_T \in H^s(T)\ \forall T \in \mathcal{T}\}$$

and are equipped with the norm $\|v\|_{H^s(\mathcal{T})}^2 := \sum_{T \in \mathcal{T}} \|v|_T\|_{H^s(T)}^2$. We will also need the corresponding semi-norm $|\cdot|_{H^s(\mathcal{T})}$. The spaces $H_h^s(\mathcal{T})$ are the spaces $H^s(\mathcal{T})$ with norm $\|v\|_{H_h^s(\mathcal{T})}^2 := \sum_{T \in \mathcal{T}} \|v|_T\|_{H_h^s(T)}^2$. For $s \in [0, 1/2)$, the space $H^s(\mathcal{T})$ is densely embedded in $L^2(\Gamma)$ such that the $L^2$-inner product $\langle \cdot, \cdot \rangle_\Gamma$ can be extended to serve as duality between $H^s(\mathcal{T})$ and $\widetilde{H}^{-s}(\mathcal{T})$. As previously mentioned, spaces with bold symbols denote vector-valued versions of scalar-valued spaces, e.g., $\mathbf{H}_h^s(\mathcal{T})$. Their elements are also denoted by bold symbols.

For a partition $\{\omega_j\}$ of $\omega \subset \Gamma$ into non-overlapping Lipschitz sub-domains (sub-surfaces), and $s \in [0,1]$, we will make frequent use of the estimate

$$\sum_j \|u|_{\omega_j}\|_{H^s(\omega_j)}^2 \leq \|u\|_{H^s(\omega)}^2 \quad \forall u \in H^s(\omega) \tag{6}$$

(which is immediate by the definition of the norms) and the following bound from [1, Theorem 4.1] (see [28], [21] for estimates in different, equivalent norms),

$$\|u\|_{\widetilde{H}^s(\omega)}^2 \lesssim \sum_j \|u|_{\omega_j}\|_{\widetilde{H}^s(\omega_j)}^2 \quad \forall u \in \widetilde{H}^s(\omega) \text{ with } u|_{\omega_j} \in \widetilde{H}^s(\omega_j)\ \forall j. \tag{7}$$

We recall some properties of the integral operators, cf. [14, 32, 34, 35].

**Proposition 2.** *For $s \in [-1/2, 1/2]$ the following mappings are continuous,*

$$\mathcal{W}: \widetilde{H}^{1/2+s}(\Gamma) \to H^{-1/2+s}(\Gamma),$$
$$\mathcal{V}: \widetilde{\mathbf{H}}^{-1/2+s}(\Gamma) \to \mathbf{H}^{1/2+s}(\Gamma). \tag{8}$$

*Furthermore, both operators are elliptic,*

$$\langle \mathcal{W}v, v \rangle_\Gamma + m_\Gamma(v,v) \gtrsim \|v\|_{\widetilde{H}^{1/2}(\Gamma)}^2 \quad \forall v \in \widetilde{H}^{1/2}(\Gamma), \tag{9}$$
$$\langle \mathcal{V}\mathbf{v}, \mathbf{v} \rangle_\Gamma \gtrsim \|\mathbf{v}\|_{\widetilde{\mathbf{H}}^{-1/2}(\Gamma)}^2 \quad \forall \mathbf{v} \in \widetilde{\mathbf{H}}^{-1/2}(\Gamma)$$

*and are related by*

$$\langle \mathcal{W}v, w \rangle_\Gamma = \langle \mathcal{V}\mathbf{curl}\,v, \mathbf{curl}\,w \rangle_\Gamma \quad \forall v, w \in \widetilde{H}^{1/2}(\Gamma). \tag{10}$$

We need the following regularity parameter (or norm of the inverse of $\mathcal{W}$):

$$C_\Gamma(s) := \sup\left\{\frac{\|\phi\|_{\widetilde{H}^{1/2+s}(\Gamma)}}{\|f\|_{H^{-1/2+s}(\Gamma)}};\ f \in H_0^{-1/2+s}(\Gamma) \setminus \{0\},\ \phi \text{ solves } (2)\right\}. \tag{11}$$

It is well known that $C_\Gamma(\cdot)$ is bounded on $[0, 1/2]$ if $\Gamma$ is a closed surface. Otherwise, $C_\Gamma(\cdot)$ is unbounded at $1/2$, cf., e.g., [38]. Therefore, in the following we restrict certain regularity parameters differently, depending on whether $\Gamma$ is closed or open. Furthermore, throughout this paper, we are interested in selecting large parameters $s$ in $C_\Gamma(s)$. Some technical details to be analyzed get more complicated when $s$ tends to 0. Therefore, subsequent stability estimate will not be uniform in $s$. To handle the different cases we define the interval

$$I_\Gamma := \begin{cases}(0, 1/2) & \text{if } \Gamma \text{ is an open surface,} \\ (0, 1/2] & \text{if } \Gamma \text{ is closed,}\end{cases}$$

and recall that

$$C_\Gamma(s) < \infty \quad \forall s \in I_\Gamma. \tag{12}$$

## 3.2 Technical results

According to [25, Cor. 1.4.4.5], for a fixed Lipschitz domain $\omega$, $H^s(\omega) = \widetilde{H}^s(\omega)$ with equivalent norms for $|s| < 1/2$. From [28, Lemma 5] we cite the following more explicit result.

**Lemma 3.** *Let $\omega \subset \mathbb{R}^2$ be a Lipschitz domain. Then there holds*

$$\|v\|_{\widetilde{H}^s(\omega)} \leq C(1/2 - |s|)^{-1}\|v\|_{H^s(\omega)} \quad \forall v \in H^s(\omega),\ \forall s \in (-1/2, 1/2)$$

*with a constant $C > 0$ that depends only on $\omega$.*

The next lemma presents scaling properties of some fractional-order Sobolev (semi-) norms.

**Lemma 4.** *Suppose that $F : \mathbb{R}^2 \to \mathbb{R}^3$ is an affine map with $\omega = F(\widehat{\omega})$ an open surface and with $\widehat{\omega}$ being a reference element. Then, for $s \in [0, 1]$,*

$$\begin{aligned}\|v\|_{\widetilde{H}^s(\omega)}^2 &\simeq h_\omega^{2-2s}\|v \circ F\|_{\widetilde{H}^s(\widehat{\omega})}^2, \\ \|v\|_{H_h^s(\omega)}^2 &\simeq h_\omega^{2-2s}\|v \circ F\|_{H^s(\widehat{\omega})}^2\end{aligned} \tag{13}$$

*and, for $s \in (0, 1)$,*

$$|v|_{H^s(\omega)}^2 \simeq h_\omega^{2-2s}|v \circ F|_{H^s(\widehat{\omega})}^2. \tag{14}$$

*Furthermore, $H_h^s(\omega) = \widetilde{H}^s(\omega)$ in a topological sense for $0 \leq s < 1/2$, and*

$$\|v\|_{\widetilde{H}^s(\omega)} \lesssim (1/2 - s)^{-1}\|v\|_{H_h^s(\omega)}. \tag{15}$$

*If $0 \leq s < t \leq 1$, then*

$$\|v\|_{H_h^s(\omega)} \lesssim h_\omega^{t-s}\|v\|_{H_h^t(\omega)}. \tag{16}$$

*All the hidden constants depend only on $\gamma_\omega$.*

*Proof.* The ubiquitous scaling results (13), (14) follow readily, cf. [28] for (13), and see [29] for transformation properties of semi-norms. According to [25, Cor. 1.4.4.5] there holds $H^s(\omega) = \widetilde{H}^s(\omega)$ and hence, $H^s_h(\omega) = \widetilde{H}^s(\omega)$ by definition of $H^s_h(\omega)$. According to Lemma 3 we have

$$\|v \circ F\|_{\widetilde{H}^s(\widehat{\omega})} \lesssim (1/2 - s)^{-1} \|v \circ F\|_{H^s(\widehat{\omega})}$$

with a constant which depends only on $\widehat{\omega}$. Finally, an application of the estimates (13) shows (15). The estimate (16) is obvious on the reference element $\widehat{\omega}$ and is transferred to $\omega$ with the estimates (13). □

The next lemma presents some estimates for different Sobolev norms and establishes bounds for surface differential operators. Here, $\mathbf{curl}_{\mathcal{T}}$ is the piecewise $\mathbf{curl}$-operator, i.e., $\mathbf{curl}_{\mathcal{T}} v|_T = \mathbf{curl}\, v|_T$ for any $T \in \mathcal{T}$ and sufficiently smooth function $v$.

**Lemma 5.** *For $s \in (0, 1)$ and $T \in \mathcal{T}$ there holds*

$$\|v\|_{H^s(T)} \lesssim \|v\|_{H^s_h(T)} \lesssim h_T^{-s} \|v\|_{H^s(T)} \quad \forall v \in H^s(T). \tag{17}$$

*Furthermore, for $s \in (-1/2, 0]$, one has*

$$h^{-s} \|\boldsymbol{\tau}\|_{\widetilde{\mathbf{H}}^s(\mathcal{T})} \lesssim \|\boldsymbol{\tau}\|_{\widetilde{\mathbf{H}}^s_h(\mathcal{T})} \lesssim \|\boldsymbol{\tau}\|_{\widetilde{\mathbf{H}}^s(\mathcal{T})} \quad \forall \boldsymbol{\tau} \in \widetilde{\mathbf{H}}^s(\mathcal{T}). \tag{18}$$

*For $s \in (0, 1/2]$ there holds*

$$\|\mathrm{curl}\,\boldsymbol{\tau}\|_{H^{-1/2+s}(\Gamma)} \lesssim \|\boldsymbol{\tau}\|_{\mathbf{H}^{1/2+s}(\Gamma)} \quad \forall \boldsymbol{\tau} \in \mathbf{H}^{1/2+s}(\Gamma), \tag{19}$$

$$\|\mathbf{curl}_{\mathcal{T}} v\|_{\widetilde{\mathbf{H}}^{-1/2+s}_h(\mathcal{T})} \lesssim |v|_{H^{1/2+s}(\mathcal{T})} \quad \forall v \in H^{1/2+s}(\mathcal{T}) \tag{20}$$

*and*

$$\|\mathbf{curl}_{\mathcal{T}} v\|_{\mathbf{H}^{-1/2+s}(\mathcal{T})} \lesssim |v|_{H^{1/2+s}(\mathcal{T})} \quad \forall v \in H^{1/2+s}(\mathcal{T}). \tag{21}$$

*All hidden constants are independent of $\mathcal{T}$ for bounded shape-regularity parameter $\gamma_{\mathcal{T}}$.*

*Proof.* The estimates (17) follow readily, and by duality we then conclude (18). To prove (19), we decompose the polyhedral surface $\Gamma$ into its faces $\Gamma_j$ and use (7) and Lemma 3 to bound

$$\|\mathrm{curl}\,\boldsymbol{\tau}\|^2_{H^{-1/2+s}(\Gamma)} \lesssim \sum_j \|\mathrm{curl}\,\boldsymbol{\tau}\|^2_{\widetilde{H}^{-1/2+s}(\Gamma_j)} \leq \sum_j c_j s^{-2} \|\mathrm{curl}\,\boldsymbol{\tau}\|^2_{H^{-1/2+s}(\Gamma_j)} \tag{22}$$

with constants $c_j$ depending on $\Gamma_j$. Now, for any face $\Gamma_j$, $\mathbf{curl} : \widetilde{H}^{1/2}(\Gamma_j) \to \widetilde{\mathbf{H}}^{-1/2}(\Gamma_j)$ is bounded (see [22, Lemma 2.2]) so that the operator curl, being the adjoint operator of $\mathbf{curl}$, is continuous from $\mathbf{H}^{1/2}(\Gamma_j)$ to $H^{-1/2}(\Gamma_j)$. By interpolation between, e.g., $\mathbf{H}^2(\Gamma_j)$ and $\mathbf{H}^{1/2}(\Gamma_j)$, one concludes that curl : $\mathbf{H}^{1/2+s}(\Gamma_j) \to H^{-1/2+s}(\Gamma_j)$ is continuous for $s \in (0, 1/2]$. Therefore, continuing the estimate (22) and also making use of (6), we obtain

$$\|\operatorname{curl}\boldsymbol{\tau}\|^2_{H^{-1/2+s}(\Gamma)} \lesssim \sum_j \|\boldsymbol{\tau}\|^2_{\mathbf{H}^{1/2+s}(\Gamma_j)} \leq \|\boldsymbol{\tau}\|^2_{\mathbf{H}^{1/2+s}(\Gamma)} \quad \forall \boldsymbol{\tau} \in \mathbf{H}^{1/2+s}(\Gamma),$$

with a hidden constant that depends on the geometry of $\Gamma$. This is (19). It remains to show (20). To this end it is enough to bound, for $T \in \mathcal{T}$ and $v \in H^{1/2+s}(T)$,

$$\|\mathbf{curl}\, v\|_{\widetilde{\mathbf{H}}_h^{-1/2+s}(T)} \lesssim |v|_{H^{1/2+s}(T)} \quad \forall v \in H^{1/2+s}(T), \tag{23}$$

with hidden constant depending on $\gamma_{\mathcal{T}}$. The definition by duality of the norm on the left-hand side of (23) and estimates (13), (15) yield

$$\|\mathbf{curl}\, v\|_{\widetilde{\mathbf{H}}_h^{-1/2+s}(T)} \lesssim s^{-1} h_T^{1/2-s} \|\mathbf{curl}\, (v \circ F)\|_{\mathbf{H}^{-1/2+s}(\widehat{T})}.$$

The operator $\mathbf{curl}$ continuously maps $H^1(\widehat{T})$ to $\mathbf{L}^2(\widehat{T})$ and also $H^{1/2}(\widehat{T})$ to $\mathbf{H}^{-1/2}(\widehat{T})$, cf. [22, Lemma 2.1]. By interpolation, it also maps $H^{1/2+s}(\widehat{T})$ continuously to $\mathbf{H}^{-1/2+s}(\widehat{T})$ and hence,

$$\|\mathbf{curl}\, (v \circ F)\|_{\mathbf{H}^{-1/2+s}(\widehat{T})} \lesssim |v \circ F|_{H^{1/2+s}(\widehat{T})} \lesssim h_T^{s-1/2} |v|_{H^{1/2+s}(T)}.$$

Here, we also used a quotient-space argument to switch to the semi-norm, cf. [29], and the scaling property (14). This yields (23). The proof of bound (21) is similar to that of (20). One just has to use that the norm $\|\cdot\|_{\mathbf{H}^{-1/2+s}(T)}$ is scalable of the same order as $\|\cdot\|_{\widetilde{\mathbf{H}}_h^{-1/2+s}(T)}$. □

For $s \in [0, 1/2)$, define the space $\mathbf{H}^s(\operatorname{curl},\omega) := \{\boldsymbol{\sigma} \in \mathbf{H}^s(\omega); \operatorname{curl}\boldsymbol{\sigma} \in H^s(\omega)\}$ with the norm $\|\boldsymbol{\sigma}\|^2_{\mathbf{H}^s(\operatorname{curl},\omega)} := \|\boldsymbol{\sigma}\|^2_{\mathbf{H}^s(\omega)} + \|\operatorname{curl}\boldsymbol{\sigma}\|^2_{H^s(\omega)}$. As before, we denote by $\boldsymbol{t}$ a generic unit tangential vector along the boundary of $\omega$, with mathematically positive orientation consistent with the orientation of $\omega$ as a subset of $\partial\Omega$. It is well known that the tangential trace mapping $\boldsymbol{\sigma} \mapsto \boldsymbol{\sigma}\cdot\boldsymbol{t}|_{\partial\omega}$ defined by

$$\langle \boldsymbol{\sigma}\cdot\boldsymbol{t}, v\rangle_{\partial\omega} := \langle \operatorname{curl}\boldsymbol{\sigma}, v\rangle_\omega - \langle \boldsymbol{\sigma}, \mathbf{curl}\, v\rangle_\omega \quad \forall v \in H^1(\omega) \tag{24}$$

is a surjective continuous mapping from $\mathbf{H}^0(\operatorname{curl},\omega)$ to $H^{-1/2}(\partial\omega)$, cf. [23, Thm. 2.5]. As mentioned previously, $H^s(\partial\omega)$ with $s < 0$ is the dual space of $H^{-s}(\partial\omega)$, the latter being the trace space of $H^{-s+1/2}(\omega)$. We will need continuity of the tangential trace mapping (without specific bounds) for spaces of higher regularity, as stated in the following lemma.

**Lemma 6.** *Let $\omega$ be a smooth surface piece with Lipschitz boundary. For $s \in (0, 1/2]$, the tangential trace mapping indicated by $\boldsymbol{\sigma} \mapsto \boldsymbol{\sigma}\cdot\boldsymbol{t}|_{\partial\omega}$ defines a surjective, linear and continuous mapping from $\mathbf{H}^{1/2-s}(\operatorname{curl},\omega)$ to $H^{-s}(\partial\omega)$.*

*Proof.* To see that the tangential trace mapping defined by (24) can be extended to a continuous application from $\mathbf{H}^{1/2-s}(\operatorname{curl},\omega)$ to $H^{-s}(\partial\omega)$ we bound by duality

$$\begin{aligned}|\langle \boldsymbol{\sigma}, \mathbf{curl}\, v\rangle_\omega| &\lesssim \|\boldsymbol{\sigma}\|_{\widetilde{\mathbf{H}}^{1/2-s}(\omega)} \|\mathbf{curl}\, v\|_{\mathbf{H}^{-1/2+s}(\omega)},\\ |\langle \operatorname{curl}\boldsymbol{\sigma}, v\rangle_\omega| &\lesssim \|\operatorname{curl}\boldsymbol{\sigma}\|_{\widetilde{H}^{1/2-s}(\omega)} \|v\|_{H^{-1/2+s}(\omega)}.\end{aligned} \tag{25}$$

The stated continuity then follows by noting that $\widetilde{H}^{1/2-s}(\omega) = H^{1/2-s}(\omega)$ (by Lemma 3), $\mathbf{curl} : H^{1/2+s}(\omega) \to \mathbf{H}^{-1/2+s}(\omega)$ (by interpolation between $\mathbf{curl} : H^1(\omega) \to \mathbf{L}^2(\omega)$ and $\mathbf{curl} : H^{1/2}(\omega) \to \mathbf{H}^{-1/2}(\omega)$, cf. [22, Lemma 2.1]), the continuous injection $H^{1/2+s}(\omega) \hookrightarrow H^{-1/2+s}(\omega)$, and the definition of $H^s(\partial\omega)$ as the trace of $H^{1/2+s}(\omega)$. The surjectivity of the tangential trace mapping follows by the existence of a solution to $\Delta\psi + \psi = 0$ on $\omega$ with $\Delta\psi := \operatorname{curl}\mathbf{curl}\,\psi$ and Neumann datum $\boldsymbol{\sigma}\cdot\boldsymbol{t}|_{\partial\omega}$, its continuous dependence $\|\psi\|_{H^{3/2-s}(\omega)} \lesssim \|\boldsymbol{\sigma}\cdot\boldsymbol{t}\|_{H^{-s}(\partial\omega)}$, and by defining $\boldsymbol{\sigma} = \mathbf{curl}\,\psi$. □

It remains to introduce yet another space. For $s \in (0, 1/2]$, we define

$$H^{-s}(\mathcal{S}) := \Big\{\hat\sigma \in \prod_{T\in\mathcal{T}} H^{-s}(\partial T);\ \exists\boldsymbol{\sigma} \in \mathbf{H}^{1/2-s}(\mathbf{curl},\Gamma) \text{ such that } \hat\sigma|_{\partial T} = \boldsymbol{\sigma}\cdot\boldsymbol{t}|_{\partial T}\ \forall T \in \mathcal{T}\Big\}.$$

It is equipped with the norm

$$\|\hat\sigma\|_{H^{-s}(\mathcal{S})} := \inf_{\boldsymbol{\sigma}\cdot\boldsymbol{t}|_{\mathcal{S}}=\hat\sigma} \|\boldsymbol{\sigma}\|_{\mathbf{H}^{1/2-s}(\mathbf{curl},\Gamma)}.$$

Here, $\boldsymbol{\sigma}\cdot\boldsymbol{t}|_{\mathcal{S}} = \hat\sigma$ means that $\boldsymbol{\sigma}\cdot\boldsymbol{t}|_{\partial T} = \hat\sigma|_{\partial T}$ for any $T \in \mathcal{T}$. For functions $\hat\sigma \in H^{-s}(\mathcal{S})$ with $\hat\sigma = \boldsymbol{\sigma}\cdot\boldsymbol{t}|_{\mathcal{S}}$ for $\boldsymbol{\sigma} \in \mathbf{H}^{1/2-s}(\mathbf{curl},\Gamma)$, and $v \in H^{1/2+s}(\mathcal{T})$ we will make use of the duality

$$\langle\hat\sigma, [v]\rangle_{\mathcal{S}} := \langle\boldsymbol{\sigma}\cdot\boldsymbol{t}, [v]\rangle_{\mathcal{S}} := \sum_{T\in\mathcal{T}} \langle\boldsymbol{\sigma}\cdot\boldsymbol{t}|_T, v|_T\rangle_{\partial T}.$$

We have the following integration-by-parts formula, cf. [22, Lemma 4.2].

**Lemma 7** (integration by parts). *Let $\Gamma$ be a Lipschitz surface. Then, for $\phi \in \widetilde{H}^{1/2}(\Gamma)$ and $\boldsymbol{\psi} \in \mathbf{H}^{1/2}(\Gamma)$, there holds*

$$\langle\mathbf{curl}\,\phi, \boldsymbol{\psi}\rangle_\Gamma = \langle\phi, \operatorname{curl}\boldsymbol{\psi}\rangle_\Gamma. \tag{26}$$

*Furthermore, if $s \in (0, 1/2]$, $v \in H^{1/2+s}(\mathcal{T})$, and $\boldsymbol{\sigma} \in \mathbf{H}^{1/2+\delta}(\Gamma)$ for some $0 < \delta < s$, it holds that*

$$\langle\operatorname{curl}\boldsymbol{\sigma}, v\rangle_\Gamma = \langle\boldsymbol{\sigma}, \mathbf{curl}_{\mathcal{T}} v\rangle_\Gamma + \langle\boldsymbol{\sigma}\cdot\boldsymbol{t}, [v]\rangle_{\mathcal{S}}. \tag{27}$$

*Proof.* The first statement follows immediately as curl is the adjoint of **curl**. Identity (27) is an extension of [22, Lemma 4.2]. It follows by using the Piola transformation on every element, integrating by parts on a reference element with the help of [22, Lemma 4.2], transforming back, and summing over the elements. See, e.g., [36, Eq.(6.22)–(6.25)] for the mapping properties involving differential operators (our case follows by rotation). □

We will measure jumps of functions $v \in H^{1/2+s}(\mathcal{T})$ in the dual space of $H^{-s}(\mathcal{S})$,

$$\|\,[v]\,\|_{H^s(\mathcal{S})} := \sup_{\boldsymbol{\sigma}\in\mathbf{H}^{1/2-s}(\mathbf{curl},\Gamma)} \frac{\langle\boldsymbol{\sigma}\cdot\boldsymbol{t}, [v]\rangle_{\mathcal{S}}}{\|\boldsymbol{\sigma}\|_{\mathbf{H}^{1/2-s}(\mathbf{curl},\Gamma)}}.$$

There holds the following trace theorem.

**Lemma 8.** *For $s \in (0, 1/2]$, we have*
$$\| [v] \|_{H^s(\mathcal{S})} \lesssim \|v\|_{H^{1/2+s}(\mathcal{T})} \quad \forall v \in H^{1/2+s}(\mathcal{T}).$$

*The hidden constant is independent of $\mathcal{T}$ for bounded shape-regularity parameter $\gamma_\mathcal{T}$.*

*Proof.* The proof is a simple combination of relation (27), dualities (25), the continuity (21), and (when $\Gamma$ is open) Lemma 3 to bound $\|\boldsymbol{\sigma}\|_{\widetilde{\mathbf{H}}^{1/2-s}(\Gamma)} \lesssim \|\boldsymbol{\sigma}\|_{\mathbf{H}^{1/2-s}(\Gamma)}$ and $\|\operatorname{curl} \boldsymbol{\sigma}\|_{\widetilde{H}^{1/2-s}(\Gamma)} \lesssim \|\operatorname{curl} \boldsymbol{\sigma}\|_{H^{1/2-s}(\Gamma)}$. □

## 4 Discontinuous Petrov-Galerkin method

This is the central section of the paper. We present and analyze a general ultra-weak formulation of the hypersingular integral equation (1) and, based on this formulation, propose a DPG method and prove its quasi-optimal convergence. The structure is as follows. In the first subsection we present the ultra-weak formulation and show that its bilinear form is definite. This is essential for the definition of norms by duality via the bilinear form. In Section 4.2 we analyze these norms, in the trial and the test space, and show essential mutual estimates (Lemmas 12 and 13). Section 4.3 establishes the equivalence of the ultra-weak formulation and the standard weak formulation (Theorem 14) and proves well-posedness of the former (Theorem 15). Finally, in Section 4.4 we recall the DPG method (in the now general setting) and present the main result (Theorem 17) which establishes a general Céa estimate for the DPG method.

### 4.1 General ultra-weak formulation

In the following, we assume that the right-hand side function $f$ of the hypersingular integral equation is $L^2(\Gamma)$-regular. This restriction can be relaxed to regularity of negative order in a relatively straightforward way. However, stability of the ultra-weak formulation in this case requires that $f \in \widetilde{H}^{-1/2+s}(\mathcal{T})$ (the dual space of test functions for $f$), which does not seem to be a practical condition when $s < 1/2$. We therefore restrict ourselves to $f \in L^2_0(\Gamma)$.

Now, for a given parameter $s \in (0, 1/2]$, the general ultra-weak formulation of (2) reads as follows: Find
$$(\boldsymbol{\sigma}, \phi, \hat{\sigma}) \in U^s := \mathbf{H}_h^{1/2-s}(\mathcal{T}) \times \widetilde{H}^{1/2-s}(\Gamma) \times H^{-s}(\mathcal{S})$$
such that
$$\langle \boldsymbol{\sigma}, \boldsymbol{\tau} \rangle_\Gamma - \langle \phi, \operatorname{curl} \mathcal{V}\boldsymbol{\tau} \rangle_\Gamma = 0 \tag{28a}$$
$$\langle \boldsymbol{\sigma}, \operatorname{\mathbf{curl}}_\mathcal{T} v \rangle_\Gamma + \langle \hat{\sigma}, [v] \rangle_\mathcal{S} + m_\Gamma(\phi, v) = \langle f, v \rangle_\Gamma \tag{28b}$$
for all $(\boldsymbol{\tau}, v) \in V^s$ with
$$V^s := \widetilde{\mathbf{H}}_h^{-1/2+s}(\mathcal{T}) \times H^{1/2+s}(\mathcal{T}).$$
Recall that the bilinear form $m_\Gamma(\cdot, \cdot)$ has been defined in (3). As in (4), we abbreviate (28) by
$$(\boldsymbol{\sigma}, \phi, \hat{\sigma}) \in U^s: \quad b(\boldsymbol{\sigma}, \phi, \hat{\sigma}; \boldsymbol{\tau}, v) = \langle f, v \rangle_\Gamma \quad \forall (\boldsymbol{\tau}, v) \in V^s,$$

where the bilinear form $b: U^s \times V^s \to \mathbb{R}$ is defined as

$$b(\boldsymbol{\sigma}, \phi, \hat{\sigma}; \boldsymbol{\tau}, v) := \langle \boldsymbol{\sigma}, \boldsymbol{\tau} + \mathbf{curl}_\mathcal{T} v \rangle_\Gamma - \langle \phi, \operatorname{curl} \mathcal{V} \boldsymbol{\tau} \rangle_\Gamma + \langle \hat{\sigma}, [v] \rangle_\mathcal{S} + m_\Gamma(\phi, v).$$

Again, the spaces $U^s$ and $V^s$ are Hilbert spaces with norms

$$\|(\boldsymbol{\sigma}, \phi, \hat{\sigma})\|_{U^s, \alpha, \beta, \rho} := \left( \alpha^2 \|\boldsymbol{\sigma}\|^2_{\mathbf{H}_h^{1/2-s}(\mathcal{T})} + \beta^2 \|\phi\|^2_{\widetilde{H}^{1/2-s}(\Gamma)} + \rho^2 \|\hat{\sigma}\|^2_{H^{-s}(\mathcal{S})} \right)^{1/2} \tag{29}$$

(depending on three parameters $\alpha, \beta, \rho > 0$) and

$$\|(\boldsymbol{\tau}, v)\|_{V^s} := \left( \|\boldsymbol{\tau}\|^2_{\widetilde{\mathbf{H}}_h^{-1/2+s}(\mathcal{T})} + \|v\|^2_{H^{1/2+s}(\mathcal{T})} \right)^{1/2}. \tag{30}$$

**Lemma 9.** *For $s \in (0, 1/2]$ the bilinear form $b(\cdot, \cdot)$ is definite, i.e.,*

$$b(\boldsymbol{\sigma}, \phi, \hat{\sigma}; \boldsymbol{\tau}, v) = 0 \quad \forall (\boldsymbol{\sigma}, \phi, \hat{\sigma}) \in U^s \Leftrightarrow (\boldsymbol{\tau}, v) = 0, \tag{31a}$$
$$b(\boldsymbol{\sigma}, \phi, \hat{\sigma}; \boldsymbol{\tau}, v) = 0 \quad \forall (\boldsymbol{\tau}, v) \in V^s \Leftrightarrow (\boldsymbol{\sigma}, \phi, \hat{\sigma}) = 0. \tag{31b}$$

*Proof.* Both $(\boldsymbol{\tau}, v) = 0$ and $(\boldsymbol{\sigma}, \phi, \hat{\sigma}) = 0$ yield $b(\boldsymbol{\sigma}, \phi, \hat{\sigma}; \boldsymbol{\tau}, v) = 0$, hence it remains to show the reverse implications. For (31a), choose $(\boldsymbol{\sigma}, \phi, \hat{\sigma}) = (\boldsymbol{\sigma}, 0, 0)$ with arbitrary $\boldsymbol{\sigma} \in \mathbf{H}^{1/2-s}(\Gamma)$. Then, $\langle \boldsymbol{\sigma}, \boldsymbol{\tau} + \mathbf{curl}_\mathcal{T} v \rangle_\Gamma = 0$, and due to $\mathbf{curl}_\mathcal{T} v \in \mathbf{H}^{-1/2+s}(\Gamma)$, it holds that $\boldsymbol{\tau} = -\mathbf{curl}_\mathcal{T} v$ as functionals in $\widetilde{\mathbf{H}}^{-1/2+s}(\Gamma)$. Next, we show that $v \in \widetilde{H}^{1/2+s}(\Gamma)$. To this end, choose a pair of elements $T_1, T_2 \in \mathcal{T}$ such that $\emptyset \ne e := (\overline{T}_1 \cap \overline{T}_2)^\circ$. Here, $M^\circ$ denotes the interior of a set $M$. Choose $\hat{\sigma} \in \widetilde{H}^{-s}(e)$ and extend it by zero to $\hat{\sigma}_1 \in \widetilde{H}^{-s}(\partial T_1)$ and $\hat{\sigma}_2 \in \widetilde{H}^{-s}(\partial T_2)$. Define the function $\Phi$ by

$$\Phi := \begin{cases} 0 & \text{on } \Gamma \setminus (T_1 \cup T_2) \\ \Phi_i & \text{on } T_i, i = 1, 2, \end{cases}$$

where $\Phi_i \in \mathbf{H}^{1/2-s}(\operatorname{curl}, T_i)$, $i = 1, 2$, are such that $\Phi_i \cdot \boldsymbol{t}|_{\partial T_i} = \hat{\sigma}_i$, cf. Lemma 6. Due to the continuity of the tangential components on $\partial T_i$, $i = 1, 2$, it holds that $\Phi \in \mathbf{H}^{1/2-s}(\operatorname{curl}, \Gamma)$. This shows that

$$0 = \langle \Phi \cdot \boldsymbol{t}_{\partial T_1}, v \rangle_{\partial T_1} + \langle \Phi \cdot \boldsymbol{t}_{\partial T_2}, v \rangle_{\partial T_2} = \langle \hat{\sigma}, [v] \rangle_e \quad \forall \hat{\sigma} \in \widetilde{H}^{-s}(e),$$

hence $[v] = 0 \in H^s(e)$. We conclude that $v \in H^{1/2+s}(\Gamma)$. For $\emptyset \ne e := T \cap \partial \Gamma$, the same argument yields $v = 0 \in H^s(e)$, i.e., $v \in \widetilde{H}^{1/2+s}(\Gamma)$. To summarize the two preceding arguments, we have $v \in \widetilde{H}^{1/2+s}(\Gamma)$ and $\boldsymbol{\tau} = -\mathbf{curl}\, v$, so that the choice $(\boldsymbol{\sigma}, \phi, \hat{\sigma}) = (0, v, 0)$ and the integration-by-parts formula (26) show

$$0 = -\langle v, \operatorname{curl} \mathcal{V} \boldsymbol{\tau} \rangle_\Gamma + m_\Gamma(v, v) = \langle \mathbf{curl}\, v, \mathcal{V} \mathbf{curl}\, v \rangle_\Gamma + m_\Gamma(v, v).$$

Relations (9) and (10) imply that $v = 0$. This also means that $\boldsymbol{\tau} = 0$. To see the direction "$\Rightarrow$" in (31b), choosing $(\boldsymbol{\tau}, v) = (0, 1)$ in case of closed $\Gamma$ shows that $\phi$ has vanishing integral mean.

Then, for arbitrary $w \in \widetilde{H}^{1/2+s}(\Gamma)$, we choose the test function $(\boldsymbol{\tau}, v) = (\mathbf{curl}\, w, -w)$. Taking into account the vanishing integral mean of $\phi$ in case of $\Gamma$ being closed, this yields

$$\langle \phi, \mathrm{curl}\, \mathcal{V}\mathbf{curl}\, w \rangle_\Gamma = \langle \mathcal{W}\phi, w \rangle_\Gamma = 0,$$

that is, $\mathcal{W}\phi = 0$ in $H^{-1/2-s}(\Gamma)$ and we conclude that $\phi = 0$. Now, it follows easily by duality arguments that $\boldsymbol{\sigma} = 0$ and $\hat{\sigma} = 0$. □

## 4.2 Norm equivalences

In Lemma 12, we will investigate the relation between the norm $\|\cdot\|_{V^s}$ from (30) and

$$\|(\boldsymbol{\tau}, v)\|_{V^s,\mathrm{opt},\alpha,\beta,\rho} := \sup_{(\boldsymbol{\sigma},\phi,\hat{\sigma}) \in U^s \setminus \{0\}} \frac{b(\boldsymbol{\sigma}, \phi, \hat{\sigma}; \boldsymbol{\tau}, v)}{\|(\boldsymbol{\sigma}, \phi, \hat{\sigma})\|_{U^s,\alpha,\beta,\rho}}. \tag{32}$$

By Lemma 9 this defines a norm for $s \in (0, 1/2]$. Inspection shows that

$$\|(\boldsymbol{\tau}, v)\|_{V^s,\mathrm{opt},\alpha,\beta,\rho} \simeq$$
$$\alpha^{-1} \|\boldsymbol{\tau} + \mathbf{curl}_\mathcal{T} v\|_{\widetilde{\mathbf{H}}_h^{s-1/2}(\mathcal{T})} + \beta^{-1} \|\mathrm{curl}\, \mathcal{V}\boldsymbol{\tau}\|_{H^{s-1/2}(\Gamma)} + \rho^{-1} \| [v] \|_{H^s(\mathcal{S})} + \beta^{-1} m_\Gamma(v,v)^{1/2}.$$

Recall that, by definition, $m_\Gamma(v,v)^{1/2} = |\langle v, 1 \rangle_\Gamma|$ if $\Gamma$ is closed, and $m_\Gamma(v,v) = 0$ otherwise. By duality arguments, the norm equivalence $\|\cdot\|_{V^s} \simeq \|\cdot\|_{V^s,\mathrm{opt},\alpha,\beta,\rho}$ we aim at implies equivalence in $U^s$ of the norm $\|\cdot\|_{U^s,\alpha,\beta,\rho}$ from (29) and the so-called energy norm

$$\|(\boldsymbol{\sigma}, \phi, \hat{\sigma})\|_{U^s} := \sup_{(\boldsymbol{\tau},v) \in V^s \setminus \{0\}} \frac{b(\boldsymbol{\sigma}, \phi, \hat{\sigma}; \boldsymbol{\tau}, v)}{\|(\boldsymbol{\tau}, v)\|_{V^s}}, \tag{33}$$

which is a norm for $s \in (0, 1/2]$ due to Lemma 9. Proving the equivalence of norms in the test space requires studying the stability of the adjoint problem. This will be done in the next two lemmas.

**Lemma 10.** *Given $s \in I_\Gamma$, $\mathbf{g}_1 \in \widetilde{\mathbf{H}}_h^{-1/2+s}(\mathcal{T})$ and $g_2 \in H_0^{-1/2+s}(\Gamma)$, there is a unique weak solution $(\boldsymbol{\tau}_c, v_c) \in \widetilde{\mathbf{H}}_h^{-1/2+s}(\mathcal{T}) \times \widetilde{H}_0^{1/2+s}(\Gamma)$ of the problem*

$$\begin{aligned}\boldsymbol{\tau}_c + \mathbf{curl}\, v_c &= \mathbf{g}_1, \\ \mathrm{curl}\, \mathcal{V}\boldsymbol{\tau}_c &= g_2.\end{aligned} \tag{34}$$

*In addition,*

$$\|v_c\|_{\widetilde{H}^{1/2+s}(\Gamma)} + \|\boldsymbol{\tau}_c\|_{\widetilde{\mathbf{H}}_h^{-1/2+s}(\mathcal{T})} \lesssim C_\Gamma(s) \left( h^{-1/2+s} \|\mathbf{g}_1\|_{\widetilde{\mathbf{H}}_h^{-1/2+s}(\mathcal{T})} + \|g_2\|_{H^{-1/2+s}(\Gamma)} \right)$$

*with $C_\Gamma(s)$ being the regularity parameter defined in (11). The hidden constant in the estimate above is independent of $\mathcal{T}$ as long as the shape-regularity parameter $\gamma_\mathcal{T}$ is bounded.*

*Proof.* Since $\widetilde{\mathbf{H}}_h^{-1/2+s}(\mathcal{T}) \subset \widetilde{\mathbf{H}}^{-1/2+s}(\Gamma)$ by (7) there holds $\mathbf{g}_1 \in \widetilde{\mathbf{H}}^{-1/2+s}(\Gamma)$ and hence, $\operatorname{curl} \mathcal{V} \mathbf{g}_1 - g_2 \in H_0^{-1/2+s}(\Gamma)$. Therefore, combining the equations (34) and using (10), we find that $v_c \in \widetilde{H}_0^{1/2}(\Gamma)$ solves the hypersingular integral equation

$$\langle \operatorname{curl} \mathcal{V} \mathbf{curl}\, v_c, \psi \rangle = \langle \mathcal{W} v_c, \psi \rangle = \langle \operatorname{curl} \mathcal{V} \mathbf{g}_1 - g_2, \psi \rangle \qquad \forall \psi \in \widetilde{H}^{1/2}(\Gamma).$$

By (9), $v_c$ is unique. Regularity theory for this problem, cf. (11), (12), shows that $v_c \in \widetilde{H}^{1/2+s}(\Gamma)$ and

$$\begin{aligned}
\|v_c\|_{\widetilde{H}^{1/2+s}(\Gamma)} &\lesssim C_\Gamma(s) \|\operatorname{curl} \mathcal{V} \mathbf{g}_1 - g_2\|_{H^{-1/2+s}(\Gamma)} \\
&\lesssim C_\Gamma(s) \Big(h^{-1/2+s} \|\mathbf{g}_1\|_{\widetilde{\mathbf{H}}_h^{-1/2+s}(\mathcal{T})} + \|g_2\|_{H^{-1/2+s}(\Gamma)}\Big).
\end{aligned} \qquad (35)$$

Here, the last estimate follows from the triangle inequality, the mapping properties (19), (8) of curl and $\mathcal{V}$, respectively, and the first bound in (18). Due to estimates (20) and (6) there holds

$$\|\mathbf{curl}\, v_c\|_{\widetilde{\mathbf{H}}_h^{-1/2+s}(\mathcal{T})} \lesssim \|v_c\|_{H^{1/2+s}(\mathcal{T})} \leq \|v_c\|_{\widetilde{H}^{1/2+s}(\Gamma)}. \qquad (36)$$

With $\boldsymbol{\tau}_c := \mathbf{g}_1 - \mathbf{curl}\, v_c$, the pair $(\boldsymbol{\tau}_c, v_c)$ is the unique weak solution of (34), and

$$\begin{aligned}
\|\boldsymbol{\tau}_c\|_{\widetilde{\mathbf{H}}_h^{-1/2+s}(\mathcal{T})} &\lesssim \|\mathbf{g}_1\|_{\widetilde{\mathbf{H}}_h^{-1/2+s}(\mathcal{T})} + \|\mathbf{curl}\, v_c\|_{\widetilde{\mathbf{H}}_h^{-1/2+s}(\mathcal{T})} \\
&\lesssim C_\Gamma(s)\Big(h^{-1/2+s} \|\mathbf{g}_1\|_{\widetilde{\mathbf{H}}_h^{-1/2+s}(\mathcal{T})} + \|g_2\|_{H^{-1/2+s}(\Gamma)}\Big).
\end{aligned} \qquad (37)$$

Here, the last estimate follows from estimates (35), (36), and by using that $1 \lesssim h^{-1/2+s}$. Combining (35) and (37) concludes the proof. $\square$

**Lemma 11.** *For $s \in (0, 1/2]$ let $(\boldsymbol{\tau}, v) \in \widetilde{\mathbf{H}}_h^{-1/2+s}(\mathcal{T}) \times H_0^{1/2+s}(\mathcal{T})$ solve*

$$\begin{aligned}
\boldsymbol{\tau} + \mathbf{curl}_\mathcal{T} v &= 0, \\
\operatorname{curl} \mathcal{V} \boldsymbol{\tau} &= 0.
\end{aligned} \qquad (38)$$

*There is a number $C_{\mathrm{hom}}(\mathcal{T}, s) > 0$, which depends only on $\mathcal{T}$ and $s$, such that*

$$\|v\|_{H^{1/2+s}(\mathcal{T})} + \|\boldsymbol{\tau}\|_{\widetilde{\mathbf{H}}_h^{-1/2+s}(\mathcal{T})} \leq C_{\mathrm{hom}}(\mathcal{T}, s) \|[v]\|_{H^s(\mathcal{S})}.$$

*Proof.* Since $\boldsymbol{\tau} = -\mathbf{curl}_\mathcal{T} v$ and $\mathbf{curl}_\mathcal{T} : H^{1/2+s}(\mathcal{T}) \to \widetilde{\mathbf{H}}^{-1/2+s}(\mathcal{T})$ (with bound depending on $\mathcal{T}$ and $s$), it suffices to show that

$$\|v\|_{H^{1/2+s}(\mathcal{T})} \leq C_{\mathrm{hom}}(\mathcal{T}, s) \|[v]\|_{H^s(\mathcal{S})}, \qquad (39)$$

where for simplicity we use the same name for the constant. The proof is indirect. The $v$-components of solutions $(\boldsymbol{\tau}, v)$ to (38) satisfy $\mathcal{W}_\mathcal{T} v := \operatorname{curl} \mathcal{V} \mathbf{curl}_\mathcal{T} v = 0$. The operator $\mathcal{W}_\mathcal{T} : H^{1/2+s}(\mathcal{T}) \to H^{-1/2+s}(\Gamma)$ is continuous since $\mathbf{curl}_\mathcal{T} : H^{1/2+s}(\mathcal{T}) \to \mathbf{H}^{-1/2+s}(\mathcal{T}) =$ $\widetilde{\mathbf{H}}^{-1/2+s}(\mathcal{T}) \hookrightarrow \mathbf{H}^{-1/2+s}(\Gamma) = \widetilde{\mathbf{H}}^{-1/2+s}(\Gamma)$ (by (21), Lemma 3, (7), and again Lemma 3), $\mathcal{V}:\widetilde{\mathbf{H}}^{-1/2+s}(\Gamma) \to \mathbf{H}^{1/2+s}(\Gamma)$ by (8), and $\mathrm{curl}: \mathbf{H}^{1/2+s}(\Gamma) \to H^{-1/2+s}(\Gamma)$ according to (19). Hence the subspace of functions $v \in H_0^{1/2+s}(\mathcal{T})$ with $\mathcal{W}_\mathcal{T} v = 0$ is closed.

Suppose that (39) does not hold. Then there is a sequence $(v_j) \subset H_0^{1/2+s}(\mathcal{T})$ whose elements satisfy

$$\mathcal{W}_\mathcal{T} v_j = 0, \quad \|v_j\|_{H^{1/2+s}(\mathcal{T})} = 1 \quad \forall j, \quad \text{and} \quad \|\,[v_j]\,\|_{H^s(\mathcal{S})} \to 0 \quad (j \to \infty).$$

By Rellich's imbedding theorem (applied element-wise) there exists a subsequence, again denoted by $(v_j)$, that converges to an element $v \in H_0^{1/2+s/2}(\mathcal{T})$. By continuity, and since $\|\,[v_j]\,\|_{H^{s/2}(\mathcal{S})} \leq \|\,[v_j]\,\|_{H^s(\mathcal{S})} \to 0$, this element satisfies $\|\,[v]\,\|_{H^{s/2}(\mathcal{S})} = 0$, so that $v \in \widetilde{H}_0^{1/2+s/2}(\Gamma)$. Furthermore, by closedness, $\mathcal{W}_\mathcal{T} v = \mathcal{W} v = 0$. The ellipticity (9) of $\mathcal{W}$ on $\widetilde{H}_0^{1/2}(\Gamma)$ implies $v = 0$. This is a contradiction to $\|v_j\|_{H^{1/2+s}(\mathcal{T})} = 1 \; \forall j$. □

Now we are ready to prove a central norm equivalence in the test space $V^s$.

**Lemma 12.** *For $s \in (0, 1/2]$ there holds*

$$\|(\boldsymbol{\tau}, v)\|_{V^s,\mathrm{opt},\alpha,\beta,\rho} \lesssim \|(\boldsymbol{\tau}, v)\|_{V^s} \quad \forall (\boldsymbol{\tau}, v) \in V^s \tag{40}$$

*with $\alpha = 1$, $\beta = h^{-1/2+s}$, $\rho = 1$, and, for $s \in I_\Gamma$, we have*

$$\|(\boldsymbol{\tau}, v)\|_{V^s} \lesssim \|(\boldsymbol{\tau}, v)\|_{V^s,\mathrm{opt},\alpha,\beta,\rho} \quad \forall (\boldsymbol{\tau}, v) \in V^s \tag{41}$$

*with $\alpha = C_\Gamma(s)^{-1} h^{1/2-s}$, $\beta = C_\Gamma(s)^{-1}$, and $\rho = C_{\mathrm{hom}}(\mathcal{T}, s)^{-1}$. Here, $C_{\mathrm{hom}}(\mathcal{T}, s)$ is the number from Lemma 11 and $C_\Gamma(s)$ has been defined in (11). The hidden constants in the estimates above are independent of $\mathcal{T}$ as long as the shape-regularity parameter $\gamma_\mathcal{T}$ is bounded.*

*Proof.* Estimate (20) shows that

$$\|\mathbf{curl}_\mathcal{T} v\|_{\widetilde{\mathbf{H}}_h^{-1/2+s}(\mathcal{T})} \lesssim \|v\|_{H^{1/2+s}(\mathcal{T})} \quad \forall v \in H^{1/2+s}(\mathcal{T}),$$

that is,

$$\|\boldsymbol{\tau} + \mathbf{curl}_\mathcal{T} v\|_{\widetilde{\mathbf{H}}_h^{s-1/2}(\mathcal{T})} \lesssim \|\boldsymbol{\tau}\|_{\widetilde{\mathbf{H}}_h^{s-1/2}(\mathcal{T})} + \|v\|_{H^{1/2+s}(\mathcal{T})} \quad \forall (\boldsymbol{\tau}, v) \in V^s.$$

Moreover, by the continuity of curl (19) and $\mathcal{V}$ (8), using (18) and (7),

$$\|\mathrm{curl}\,\mathcal{V}\boldsymbol{\tau}\|_{H^{-1/2+s}(\Gamma)} \lesssim \|\mathcal{V}\boldsymbol{\tau}\|_{H^{1/2+s}(\Gamma)} \lesssim \|\boldsymbol{\tau}\|_{\widetilde{\mathbf{H}}^{-1/2+s}(\Gamma)} \lesssim h^{-1/2+s}\|\boldsymbol{\tau}\|_{\widetilde{\mathbf{H}}_h^{-1/2+s}(\mathcal{T})}$$

for any $\boldsymbol{\tau} \in \widetilde{\mathbf{H}}_h^{-1/2+s}(\mathcal{T})$ and $s \in (0, 1/2]$. Lemma 8 and, in the case of a closed surface, the bound $m_\Gamma(v, v)^{1/2} \leq \|v\|_{L^2(\Gamma)}\|1\|_{L^2(\Gamma)} \lesssim \|v\|_{H^{1/2+s}(\mathcal{T})} \lesssim h^{-1/2+s}\|v\|_{H^{1/2+s}(\mathcal{T})}$ conclude the proof of (40). To prove (41) we have to show that

$$\begin{aligned}\|\boldsymbol{\tau}\|_{\widetilde{\mathbf{H}}_h^{-1/2+s}(\mathcal{T})} + \|v\|_{H^{1/2+s}(\mathcal{T})} &\lesssim C_\Gamma(s) h^{-1/2+s}\|\boldsymbol{\tau} + \mathbf{curl}_\mathcal{T} v\|_{\widetilde{\mathbf{H}}_h^{s-1/2}(\mathcal{T})} \\ &+ C_\Gamma(s)\|\mathrm{curl}\,\mathcal{V}\boldsymbol{\tau}\|_{H^{s-1/2}(\Gamma)} + C_{\mathrm{hom}}(\mathcal{T}, s)\|\,[v]\,\|_{H^s(\mathcal{S})} + C_\Gamma(s) m_\Gamma(v, v)^{1/2}\end{aligned} \tag{42}$$

for any $(\boldsymbol{\tau}, v) \in V^s$ and $s \in I_\Gamma$. To this end, for given $(\boldsymbol{\tau}, v) \in V^s$ define $\mathbf{g}_1 := \boldsymbol{\tau} + \mathbf{curl}_\mathcal{T} v$ and $g_2 := \mathrm{curl}\, \mathcal{V}\boldsymbol{\tau}$, and denote by $(\boldsymbol{\tau}_c, v_c)$ the solution of (34) established by Lemma 10.

**If $\Gamma$ is open:** Then $(\boldsymbol{\tau}_0, v_0) := (\boldsymbol{\tau}, v) - (\boldsymbol{\tau}_c, v_c)$ solves (38). The triangle inequality and stability estimate by Lemma 10 combined with $\|v_c\|_{H^{1/2+s}(\mathcal{T})} \lesssim \|v_c\|_{\widetilde{H}^{1/2+s}(\Gamma)}$ by (6), and Lemma 11 show (42).

**If $\Gamma$ is closed:** Let $c_v := |\Gamma|^{-1} \langle v, 1 \rangle_\Gamma$ be the average of $v$. It satisfies

$$\|c_v\|_{H^{1/2+s}(\mathcal{T})} = \|c_v\|_{L^2(\Gamma)} \lesssim |c_v| \lesssim m_\Gamma(c_v, c_v)^{1/2} = m_\Gamma(v, v)^{1/2}.$$

Then $(\boldsymbol{\tau}_0, v_0) := (\boldsymbol{\tau}, v) - (\boldsymbol{\tau}_c, v_c) - (0, c_v)$ solves (38), and we proceed as before by taking into account the established bound for $\|c_v\|_{H^{1/2+s}(\mathcal{T})}$ and the fact that $C_\Gamma(s)$ is bounded from below by a positive constant. □

The norm equivalence from Lemma 12 implies a corresponding equivalence in $U^s$.

**Lemma 13.** *For $s \in (0, 1/2]$ there holds*

$$\|(\boldsymbol{\sigma}, \phi, \hat{\sigma})\|_{U^s} \lesssim \|\boldsymbol{\sigma}\|_{\mathbf{H}_h^{1/2-s}(\mathcal{T})} + h^{-1/2+s}\|\phi\|_{\widetilde{H}^{1/2-s}(\Gamma)} + \|\hat{\sigma}\|_{H^{-s}(\mathcal{S})}$$

*and, for $s \in I_\Gamma$, we have*

$$C_\Gamma(s)^{-1} h^{1/2-s} \|\boldsymbol{\sigma}\|_{\mathbf{H}_h^{1/2-s}(\mathcal{T})} + C_\Gamma(s)^{-1} \|\phi\|_{\widetilde{H}^{1/2-s}(\Gamma)} + C_{\mathrm{hom}}(\mathcal{T}, s)^{-1} \|\hat{\sigma}\|_{H^{-s}(\mathcal{S})} \lesssim \|(\boldsymbol{\sigma}, \phi, \hat{\sigma})\|_{U^s}$$

*for any $(\boldsymbol{\sigma}, \phi, \hat{\sigma}) \in U^s$. Here, $C_{\mathrm{hom}}(\mathcal{T}, s)$ is the number from Lemma 11 and $C_\Gamma(s)$ has been defined in (11). The hidden constants in the estimates above are independent of $\mathcal{T}$ as long as the shape-regularity parameter $\gamma_\mathcal{T}$ is bounded.*

*Proof.* Recall the definitions (29), (30) of the norms $\|\cdot\|_{U^s, \alpha, \beta, \rho}$ and $\|\cdot\|_{V^s}$, and the definitions (32), (33) of $\|\cdot\|_{V^s, \mathrm{opt}, \alpha, \beta, \rho}$ and $\|\cdot\|_{U^s}$ by duality via the bilinear form $b(\cdot, \cdot)$. According to [39, Proposition 2.1], the initial norm $\|\cdot\|_{U^s, \alpha, \beta, \rho}$ can be recovered by the iterated duality

$$\|(\boldsymbol{\sigma}, \phi, \hat{\sigma})\|_{U^s, \alpha, \beta, \rho} = \sup_{(\boldsymbol{\tau}, v) \in V^s \setminus \{0\}} \frac{b(\boldsymbol{\sigma}, \phi, \hat{\sigma}; \boldsymbol{\tau}, v)}{\|(\boldsymbol{\tau}, v)\|_{V^s, \mathrm{opt}, \alpha, \beta, \rho}}. \tag{43}$$

Due to the definition (33) of $\|\cdot\|_{U^s}$ and the duality (43), an estimate of the type $\|(\boldsymbol{\sigma}, \phi, \hat{\sigma})\|_{U^s} \leq C\|(\boldsymbol{\sigma}, \phi, \hat{\sigma})\|_{U^s, \alpha, \beta, \rho}$ is equivalent to the estimate $\|\cdot\|_{V^s, \mathrm{opt}, \alpha, \beta, \rho} \leq C\|\cdot\|_{V^s}$, and correspondingly for the inverse inequality. Therefore, the first bound of the lemma follows from (40), and the second from (41). □

### 4.3 Solvability of the ultra-weak formulation

In order to show unique solvability of the problem (28) independent of $s$, we will make use of the fact that a weak solution of the hypersingular integral equations is also an ultra-weak solution for all $s \in (0, 1/2]$.

**Theorem 14** (weak solution is general ultra-weak solution). *Let $f \in L_0^2(\Gamma)$ and denote by $\phi \in \widetilde{H}_0^{1/2}(\Gamma)$ the unique solution of the weak formulation (2). Define $\boldsymbol{\sigma} := \mathcal{V}\mathbf{curl}\,\phi$ and $\hat{\sigma} := \prod_{T\in\mathcal{T}} \boldsymbol{\sigma}\cdot\boldsymbol{t}|_{\partial T}$ (for the definition of $\boldsymbol{\sigma}\cdot\boldsymbol{t}|_{\partial T}$ cf. (24)). Then, for every $s \in (0,1/2]$, we have $(\boldsymbol{\sigma},\phi,\hat{\sigma}) \in U^s$, and $(\boldsymbol{\sigma},\phi,\hat{\sigma})$ fulfills (28).*

*Proof.* First we note that the weak solution $\phi$ of (2) is an element of $\widetilde{H}^{1/2+s}(\Gamma)$ for any $s \in I_\Gamma$. This implies that its restriction onto a face $\Gamma_j$ satisfies $\phi|_{\Gamma_j} \in H^{1/2+s}(\Gamma_j)$ so that $\mathbf{curl}\,\phi|_{\Gamma_j} \in \mathbf{H}^{-1/2+s}(\Gamma_j)$ by (21). Now, since $\mathbf{H}^{-1/2+s}(\Gamma_j) = \widetilde{\mathbf{H}}^{-1/2+s}(\Gamma_j)$ (see Lemma 3), and since $\prod_j \widetilde{\mathbf{H}}^{-1/2+s}(\Gamma_j) \subset \widetilde{\mathbf{H}}^{-1/2+s}(\Gamma) = \mathbf{H}^{-1/2+s}(\Gamma)$ by (7) and again Lemma 3, we conclude that $\mathbf{curl}\,\phi \in \mathbf{H}^{-1/2+s}(\Gamma) = \widetilde{\mathbf{H}}^{-1/2+s}(\Gamma)$. It follows from the mapping properties (8) of $\mathcal{V}$ that $\boldsymbol{\sigma} \in \mathbf{H}^{1/2+s}(\Gamma)$ and hence also $\boldsymbol{\sigma} \in \mathbf{H}_h^{1/2}(\mathcal{T})$. The symmetry of $\mathcal{V}$ and the integration-by-parts formula (26) show

$$\langle \boldsymbol{\sigma}, \boldsymbol{\tau}\rangle_\Gamma = \langle \phi, \mathrm{curl}\,\mathcal{V}\boldsymbol{\tau}\rangle_\Gamma \quad \forall \boldsymbol{\tau} \in \widetilde{\mathbf{H}}^{-1/2+s}(\Gamma).$$

This implies (28a) since $\widetilde{\mathbf{H}}_h^{-1/2+s}(\mathcal{T}) \subset \mathbf{H}^{-1/2+s}(\Gamma) = \widetilde{\mathbf{H}}^{-1/2+s}(\Gamma)$, cf. Lemma 3 and (7). To show the remaining parts note that

$$\langle f, v\rangle = \langle \mathrm{curl}\,\boldsymbol{\sigma}, v\rangle_\Gamma = \langle \boldsymbol{\sigma}, \mathbf{curl}_\mathcal{T} v\rangle_\Gamma + \langle \boldsymbol{\sigma}\cdot\boldsymbol{t}, [v]\rangle_\mathcal{S} \quad \forall v \in H^{1/2+s}(\mathcal{T}). \tag{44}$$

The first equality in (44) is due to the fact that $\mathrm{curl}\,\boldsymbol{\sigma} = f \in L^2(\Gamma)$ and $v \in H^{1/2+s}(\mathcal{T}) \subset L^2(\Gamma)$. We have seen before that $\boldsymbol{\sigma} \in \mathbf{H}^{1/2+s}(\Gamma)$. This allows for an application of the integration-by-parts formula (27), which shows the second identity in (44). □

Lemma 9 together with Babuška-Brezzi theory shows that the ultra-weak formulation (28) is uniquely solvable.

**Theorem 15** (solvability of general ultra-weak formulation). *Let $s \in I_\Gamma$ and $f \in L_0^2(\Gamma)$ be given. The ultra-weak formulation (28) has a unique solution $(\boldsymbol{\sigma},\phi,\hat{\sigma}) \in U^s$ which is independent of $s$, and*

$$\|(\boldsymbol{\sigma},\phi,\hat{\sigma})\|_{U^s,\alpha,\beta,\rho} \lesssim \|f\|_{L^2(\Gamma)}$$

*with $\alpha = C_\Gamma(s)^{-1} h^{1/2-s}$, $\beta = C_\Gamma(s)^{-1}$, and $\rho = C_{\mathrm{hom}}(\mathcal{T},s)^{-1}$. Here, $C_{\mathrm{hom}}(\mathcal{T},s)$ is the number from Lemma 11 and $C_\Gamma(s)$ has been defined in (11). The hidden constant in the estimate is independent of $f$ and $\mathcal{T}$ as long as the shape-regularity parameter $\gamma_\mathcal{T}$ is bounded.*

*Proof.* Since $\|\cdot\|_{U^s,\alpha,\beta,\rho}$ and $\|\cdot\|_{V^s,\mathrm{opt},\alpha,\beta,\rho}$ are dual norms with respect to the bilinear form $b(\cdot,\cdot)$ (cf. (43)) it follows that

$$b(\cdot,\cdot) : (U^s, \|\cdot\|_{U^s,\alpha,\beta,\rho}) \times (V^s, \|\cdot\|_{V^s,\mathrm{opt},\alpha,\beta,\rho}) \to \mathbb{R} \quad \text{is bounded,}$$

$$\sup_{(\boldsymbol{\tau},v)\in V^s\setminus\{0\}} \frac{b(\boldsymbol{\sigma},\phi,\hat{\sigma};\boldsymbol{\tau},v)}{\|(\boldsymbol{\tau},v)\|_{V^s,\mathrm{opt},\alpha,\beta,\rho}} \geq \|(\boldsymbol{\sigma},\phi,\hat{\sigma})\|_{U^s,\alpha,\beta,\rho} \quad \forall(\boldsymbol{\sigma},\phi,\hat{\sigma}) \in U^s.$$

In addition, by the definiteness (31a) of the bilinear form $b(\cdot,\cdot)$ and since $\|\cdot\|_{U^s,\alpha,\beta,\rho}$ is a norm in $U^s$, we conclude that

$$\sup_{(\boldsymbol{\sigma},\phi,\hat{\sigma})\in U^s\setminus\{0\}} \frac{b(\boldsymbol{\sigma},\phi,\hat{\sigma};\boldsymbol{\tau},v)}{\|(\boldsymbol{\sigma},\phi,\hat{\sigma})\|_{U^s,\alpha,\beta,\rho}} > 0 \quad \forall (\boldsymbol{\tau},v) \in V^s \setminus \{0\}.$$

By standard estimates we start bounding the right-hand side in (28):

$$|\langle f, v\rangle_\Gamma| \leq \|f\|_{L^2(\Gamma)} \|v\|_{L^2(\Gamma)} \leq \|f\|_{L^2(\Gamma)} \|v\|_{H^{1/2+s}(\mathcal{T})} = \|f\|_{L^2(\Gamma)} \|(0,v)\|_{V^s} \quad \forall v \in H^{1/2+s}(\mathcal{T}).$$

Lemma 12 implies that the right-hand side functional is bounded by $O(\|f\|_{L^2(\Gamma)})$ in the norm dual to $\|\cdot\|_{V^s,\mathrm{opt},\alpha,\beta,\rho}$ if we choose $\alpha = C_\Gamma(s)^{-1} h^{1/2-s}$, $\beta = C_\Gamma(s)^{-1}$, and $\rho = C_{\mathrm{hom}}(\mathcal{T},s)^{-1}$.

Babuška-Brezzi theory [6] shows that there is a unique solution $(\boldsymbol{\sigma},\phi,\hat{\sigma}) \in U^s$ of the ultra-weak formulation (28) with

$$\|(\boldsymbol{\sigma},\phi,\hat{\sigma})\|_{U^s,\alpha,\beta,\rho} \lesssim \|f\|_{L^2(\Gamma)}$$

when choosing the previously specified parameters $\alpha$, $\beta$, $\rho$. As both the weak and the ultra-weak formulation are uniquely solvable, we conclude that the solution $(\boldsymbol{\sigma},\phi,\hat{\sigma}) \in U^s$ is independent of $s$. This finishes the proof. $\square$

The assertion of Theorem 15 excludes $s = 1/2$ when $\Gamma$ is an open surface. In that case we still have stability in the energy norm.

**Corollary 16.** *Let $f \in L_0^2(\Gamma)$ and $s \in (0,1/2]$ be given. The ultra-weak formulation (28) has a unique solution $(\boldsymbol{\sigma},\phi,\hat{\sigma}) \in U^s$ which is independent of $s$, and there holds*

$$\|(\boldsymbol{\sigma},\phi,\hat{\sigma})\|_{U^s} \lesssim \|f\|_{L^2(\Gamma)}.$$

*Proof.* By definition of the energy norm (33), being dual with respect to $b(\cdot,\cdot)$ to $\|\cdot\|_{V^s}$ and being a norm by Lemma 9, and employing the estimate for the right-hand side functional from the proof of the previous theorem, we obtain

$$b(\cdot,\cdot): (U^s, \|\cdot\|_{U^s}) \times (V^s, \|\cdot\|_{V^s}) \to \mathbb{R} \quad \text{is bounded,}$$

$$\sup_{(\boldsymbol{\tau},v)\in V^s\setminus\{0\}} \frac{b(\boldsymbol{\sigma},\phi,\hat{\sigma};\boldsymbol{\tau},v)}{\|(\boldsymbol{\tau},v)\|_{V^s}} = \|(\boldsymbol{\sigma},\phi,\hat{\sigma})\|_{U^s} \quad \forall (\boldsymbol{\sigma},\phi,\hat{\sigma}) \in U^s,$$

$$\sup_{(\boldsymbol{\sigma},\phi,\hat{\sigma})\in U^s\setminus\{0\}} \frac{b(\boldsymbol{\sigma},\phi,\hat{\sigma};\boldsymbol{\tau},v)}{\|(\boldsymbol{\sigma},\phi,\hat{\sigma})\|_{U^s}} > 0 \quad \forall (\boldsymbol{\tau},v) \in V^s \setminus \{0\},$$

$$|\langle f,v\rangle_\Gamma| \leq \|f\|_{L^2(\Gamma)} \|v\|_{L^2(\Gamma)} \leq \|f\|_{L^2(\Gamma)} \|(0,v)\|_{V^s} \quad \forall v \in H^{1/2+s}(\mathcal{T}).$$

Babuška-Brezzi theory proves the assertion in the space $\overline{U}^s$ which is the completion of $U^s$ with respect to the energy norm. By Theorem 14, $(\boldsymbol{\sigma},\phi,\hat{\sigma}) \in U^s$. This finishes the proof. $\square$

### 4.4 DPG method with optimal test functions

Consider a discrete space $U_{hp} \subset U^s$ and denote by

$$\{\boldsymbol{u}_i;\ i=1,\ldots,\dim(U_{hp})\}$$

a basis of $U_{hp}$, and denote the discrete space of test functions by $V_{hp} \subset V^s$ with basis given by

$$\{\mathbf{v}_i := \Theta(\boldsymbol{u}_i);\ i=1,\ldots,\dim(U_{hp})\}, \tag{45}$$

where $\Theta : U^s \to V^s$ is the trial-to-test operator, defined by

$$\langle \Theta(\boldsymbol{u}_i), \mathbf{v}\rangle_{V^s} = b(\boldsymbol{u}_i, \mathbf{v}) \quad \forall \mathbf{v} \in V^s. \tag{46}$$

Here, $\langle \cdot, \cdot \rangle_{V^s}$ is the inner product in $V^s$ which induces the norm $\|\cdot\|_{V^s}$ in (30). The DPG method with optimal test functions as presented in [15, 16], consists in finding $(\boldsymbol{\sigma}_{hp}, \phi_{hp}, \hat{\sigma}_{hp}) \in U_{hp}$ such that

$$b(\boldsymbol{\sigma}_{hp}, \phi_{hp}, \hat{\sigma}_{hp}; \boldsymbol{\tau}, v) = \langle f, v \rangle_\Gamma \quad \forall (\boldsymbol{\tau}, v) \in V_{hp}. \tag{47}$$

A distinguishing feature of the DPG method with optimal test functions is optimal convergence in the energy norm $\|\cdot\|_{U^s}$ from (33), cf. [16, Thm. 2.2],

$$\|(\boldsymbol{\sigma} - \boldsymbol{\sigma}_{hp}, \phi - \phi_{hp}, \hat{\sigma} - \hat{\sigma}_{hp})\|_{U^s} = \inf_{(\boldsymbol{\sigma}_d, \phi_d, \hat{\sigma}_d) \in U_{hp}} \|(\boldsymbol{\sigma} - \boldsymbol{\sigma}_d, \phi - \phi_d, \hat{\sigma} - \hat{\sigma}_d)\|_{U^s}. \tag{48}$$

This best approximation property and Lemma 13 immediately imply the following quasi-optimal error estimate in standards norms.

**Theorem 17** (general Céa estimate). *Given $s \in I_\Gamma$, let $(\boldsymbol{\sigma}, \phi, \hat{\sigma}) \in U^s$ and $(\boldsymbol{\sigma}_{hp}, \phi_{hp}, \hat{\sigma}_{hp}) \in U_{hp}$ be the solutions of the ultra-weak formulation (28) and the DPG scheme (47), respectively. Then,*

$$C_\Gamma(s)^{-1} h^{1/2-s} \|\boldsymbol{\sigma} - \boldsymbol{\sigma}_{hp}\|_{\mathbf{H}_h^{1/2-s}(\mathcal{T})} + C_\Gamma(s)^{-1} \|\phi - \phi_{hp}\|_{\widetilde{H}^{1/2-s}(\Gamma)} + C_{\hom}(\mathcal{T},s)^{-1} \|\hat{\sigma} - \hat{\sigma}_{hp}\|_{H^{-s}(\mathcal{S})}$$

$$\lesssim \inf_{(\boldsymbol{\sigma}_d, \phi_d, \hat{\sigma}_d) \in U_{hp}} \left( \|\boldsymbol{\sigma} - \boldsymbol{\sigma}_d\|_{\mathbf{H}_h^{1/2-s}(\mathcal{T})} + h^{-1/2+s} \|\phi - \phi_d\|_{\widetilde{H}^{1/2-s}(\Gamma)} + \|\hat{\sigma} - \hat{\sigma}_d\|_{H^{-s}(\mathcal{S})} \right).$$

*Here, $C_{\hom}(\mathcal{T}, s)$ is the number from Lemma 11 and $C_\Gamma(s)$ has been defined in (11). The hidden constant in the estimate is independent of $\mathcal{T}$ as long as the shape-regularity parameter $\gamma_\mathcal{T}$ is bounded.*

In principle, for every $s \in I_\Gamma$, one obtains a DPG method with corresponding error estimate. Nevertheless, only the case $s = 1/2$ is practical since then, the inner product in $V^s$ reduces to $L^2$ and piecewise $\mathbf{H}^1$-bilinear forms which are easy to implement (for the calculation of test functions and for error control). On open surfaces, taking the limit $s \to 1/2$ lets $C_\Gamma(s)$ tend to infinity so that the lower bound of the error estimate from Theorem 17 is useless in this case. However, by Corollary 16, there is still a stable unique ultra-weak solution in $U^{1/2}$ and Lemma 13 provides an upper bound for the energy norm in $U^s$ for any $s \in (0, 1/2]$. Therefore, using the best approximation property (48), we can state the following error estimate that comprises the extreme case of an open surface and test space $V^{1/2}$.

**Corollary 18.** *Let $f \in L_0^2(\Gamma)$ and $s \in (0, 1/2]$ be given. Furthermore, let $(\boldsymbol{\sigma}, \phi, \hat{\sigma}) \in U^s$ and $(\boldsymbol{\sigma}_{hp}, \phi_{hp}, \hat{\sigma}_{hp}) \in U_{hp}$ be the solutions of the ultra-weak formulation (28) and the DPG scheme (47), respectively. Then,*

$$\|(\boldsymbol{\sigma} - \boldsymbol{\sigma}_{hp}, \phi - \phi_{hp}, \hat{\sigma} - \hat{\sigma}_{hp})\|_{U^s}$$
$$\lesssim \inf_{(\boldsymbol{\sigma}_d, \phi_d, \hat{\sigma}_d) \in U_{hp}} \left( \|\boldsymbol{\sigma} - \boldsymbol{\sigma}_d\|_{\mathbf{H}_h^{1/2-s}(\mathcal{T})} + h^{-1/2+s}\|\phi - \phi_d\|_{\widetilde{H}^{1/2-s}(\Gamma)} + \|\hat{\sigma} - \hat{\sigma}_d\|_{H^{-s}(\mathcal{S})} \right).$$

*The hidden constant in the estimate is independent of $\mathcal{T}$ as long as the shape-regularity parameter $\gamma_\mathcal{T}$ is bounded.*

Let us conclude the theoretical part with the following remark.

**Remark 19.** By the missing $H^1$-regularity of the solution $\phi$ to the hypersingular integral equation on open surfaces, our analysis did not lead to an estimate $\|\cdot\|_{V^{1/2}} \lesssim \|\cdot\|_{V^{1/2},\mathrm{opt},\alpha,\beta,\rho}$ for the two test norms in this case, cf. Lemma 12. For this reason, stability of the ultra-weak formulation on open surfaces, and an error estimation for the DPG method, were obtained for $s = 1/2$ only in the energy norm $\|\cdot\|_{U^{1/2}}$. This norm is weaker than the $L^2$-norm (for principal unknowns), though we did not prove that it is strictly weaker. One could use weighted test norms to obtain estimates in product norms (rather than the energy norm which mixes different components). In this way, one expects to achieve stability and error estimates in weighted $L^2$-norms (for principal unknowns). Nevertheless, in this paper we stick to standard Sobolev norms which are easier to implement than weighted Sobolev norms.

## 5 Numerical experiments

We report on four numerical experiments. The underlying surfaces are piecewise flat, and the partitions $\mathcal{T}$ consist of triangles $T$ without hanging nodes. We consider lowest order, piecewise constant functions for all components of the approximation space $U_{hp}(\mathcal{T})$. In particular, approximations for $\hat{\sigma}$ on $\mathcal{S}$ are piecewise constant on the edges of the mesh $\mathcal{T}$. Note that for all $s \in (0, 1/2]$, $U_{hp}(\mathcal{T}) \subset U^s$. The discrete test space $V_{hp}$ from (45) has finite dimension but, as discussed in the introduction, $V^s$ is infinite-dimensional so that the trial-to-test operator $\Theta$ from (46) must be approximated. This gives rise to the so-called *practical* DPG method, cf. [24]. It consists in selecting a finite-dimensional subspace $V_{h,p+r}(\mathcal{T}) \subset V^s$ and using, instead of $\Theta$, the approximated operator that one obtains by solving (46) in $V_{h,p+r}(\mathcal{T})$. As indicated by the notation, we define $V_{h,p+r}(\mathcal{T})$ as a piecewise polynomial space on the mesh $\mathcal{T}$ used for $U_{hp}$, and increase polynomial degrees. This is common procedure. In our experiments we increase degrees by 2 ($r = 2$ in the notation), which amounts to using piecewise polynomials of degree 2 for $\boldsymbol{\tau}$ (which is an $\mathbf{L}^2$-function with lowest order 0), and of degree 3 for $v$ (since the lowest order for an $H^1$-function is 1). Our choice for $s$ will always be $s = 1/2$. The inner product in the test space $V^{1/2}$ is

$$\langle \mathbf{v}, \delta_{\mathbf{v}} \rangle_{V^{1/2}} = \langle \boldsymbol{\tau}, \delta_{\boldsymbol{\tau}} \rangle_\Gamma + \langle v, \delta_v \rangle_\Gamma + \langle \mathbf{curl}_\mathcal{T} v, \mathbf{curl}_\mathcal{T} \delta_v \rangle_\Gamma$$

with $\mathbf{v} = (\boldsymbol{\tau}, v)$ and $\delta_{\mathbf{v}} = (\delta_{\boldsymbol{\tau}}, \delta_v)$. A key property of $\langle \cdot, \cdot \rangle_{V^{1/2}}$ is its locality, that is,

$$\langle \mathbf{v}, \delta_{\mathbf{v}} \rangle_{V^{1/2}} = \langle \boldsymbol{\tau}, \delta_{\boldsymbol{\tau}} \rangle_T + \langle v, \delta_v \rangle_T + \langle \mathbf{curl}\, v, \mathbf{curl}\, \delta_v \rangle_T$$

for $\mathbf{v}, \delta_{\mathbf{v}} \in V^{1/2}$ with $\mathrm{supp}(\mathbf{v}), \mathrm{supp}(\delta_{\mathbf{v}}) \subset \overline{T}$ and $T \in \mathcal{T}$. Due to this property, the calculation of the trial-to-test operator amounts to solving only local problems. In other words, the left-hand side in (46) is a block-diagonal matrix $\mathbf{V} \in \mathbb{R}^{\dim V_{h,p+2}(\mathcal{T}) \times \dim V_{h,p+2}(\mathcal{T})}$ with a fixed block size of at most 10 (owing to the local spaces for $v$) and hence is cheap to apply and invert. If $\mathbf{B} \in \mathbb{R}^{\dim V_{h,p+2}(\mathcal{T}) \times \dim U_{hp}(\mathcal{T})}$ is the discretization of the bilinear form $b$ with respect to the spaces $V_{h,p+2}(\mathcal{T})$ and $U_{hp}(\mathcal{T})$, the overall system can be written as

$$\mathbf{B}^T \mathbf{V}^{-1} \mathbf{B} u = \mathbf{V}^{-1} \mathbf{B} \mathbf{f}.$$

The numerical experiments were carried out in C++. To compute the term $\langle \phi, \mathrm{curl}\, \mathcal{V} \boldsymbol{\tau} \rangle_\Gamma$, we use integration by parts and the recurrence formulas from Maischak [33] to evaluate $\mathcal{V}$. In order to compute and store the discretization efficiently, we used the library HLib[1] to employ Hierarchical Matrices [26, 27] and the ACA algorithm [2]. The other parts are stored via standard sparse representations. The matrices $\mathbf{B}$ and $\mathbf{V}$ and the vector $\mathbf{f}$ are computed once and the system is solved by a standard iterative solver without preconditioning. In all experiments, we choose a sequence of partitions $\{\mathcal{T}_i\}_{i=0}^N$, where $\mathcal{T}_0$ is the coarsest partition and $\mathcal{T}_{i+1}$ is constructed from $\mathcal{T}_i$ by refining every triangle into 4 smaller triangles by the so-called *newest vertex bisection*, which maintains shape-regularity. In all of our experiments, we plot the error in energy norm $\|(\boldsymbol{\sigma} - \boldsymbol{\sigma}_{hp}, \phi - \phi_{hp}, \hat{\sigma} - \hat{\sigma}_{hp})\|_{U^{1/2}}$ in double logarithmic scale versus the number of triangles of the current partition. In two of the experiments we know the exact solution, such that we can also plot the errors in $L^2$-norm, $\|\boldsymbol{\sigma} - \boldsymbol{\sigma}_{hp}\|_{\mathbf{L}^2(\Gamma)}$ and $\|\phi - \phi_{hp}\|_{L^2(\Gamma)}$. Additionally, for curiosity, we plot the error $\|\hat{\sigma} - \hat{\sigma}_{hp}\|_{L^2(\mathcal{S})}$. Note that we have not proved any precise estimate for the approximation of $\hat{\sigma}$ (the number $C_{\mathrm{hom}}(\mathcal{T}, s)$ in the Céa estimate from Theorem 17 is unknown and is expected to depend on the mesh). Moreover, the natural norm for $\hat{\sigma}$ is of order $-1/2$, not 0.

Now, by Theorem 17 and standard approximation theory [36], the expected convergence order is

$$\|\boldsymbol{\sigma} - \boldsymbol{\sigma}_{hp}\|_{\mathbf{L}^2(\Gamma)} + \|\phi - \phi_{hp}\|_{L^2(\Gamma)} \lesssim h^\alpha \big(\|\boldsymbol{\sigma}\|_{\mathbf{H}^\alpha(\Gamma)} + \|\phi\|_{H^\alpha(\Gamma)} + \|f\|_{H^\alpha(\Gamma)}\big), \quad \alpha \in (0, 1].$$

Here, one uses the bound $\|\hat{\sigma} - \boldsymbol{\sigma}_d \cdot \mathbf{t}|_\mathcal{S}\|_{H^{-1/2}(\mathcal{S})} \lesssim \|\boldsymbol{\sigma} - \boldsymbol{\sigma}_d\|_{\mathbf{H}(\mathrm{curl}, \Gamma)}$ and selects an $\mathbf{H}(\mathrm{curl}, \Gamma)$-interpolant $\boldsymbol{\sigma}_d$ of $\boldsymbol{\sigma}$ (note that $\hat{\sigma} = \boldsymbol{\sigma} \cdot \mathbf{t}|_\mathcal{S}$ and $\mathrm{curl}\, \boldsymbol{\sigma} = f$.)

In all the experiments below, we observe the maximum convergence order $h = \mathcal{O}(\#\mathcal{T}^{-1/2})$. Since we plot squares of the errors, this is confirmed by the curves $\#\mathcal{T}^{-1}$.

**Experiment I.** We choose $\Gamma$ to be the surface of $\Omega$ which is a cube with side length 2, centered at the origin. The coarsest mesh $\mathcal{T}_0$ consists of 12 triangles, 2 on every side of the cube. The exact solution $\phi$ is prescribed as a piecewise affine, globally continuous function, with values at the nodes of $\mathcal{T}_0$ such that the mean value of $\phi$ vanishes. The outcome of our method with uniform mesh-refinement is shown in Figure 1. It indicates convergence with order $O(h)$. There, and in the following, $\boldsymbol{u}$ and $\boldsymbol{u}_{hp}$ stand respectively for $(\boldsymbol{\sigma}, \phi, \hat{\sigma})$ and its approximation. Note that we do not have precise information about the resulting right-hand side function $f$. It is unlikely that $f \in H^1(\Gamma)$. Though $f$ might be piecewise smooth in which case an order $O(h)$ could be proved.

---
[1] www.hlib.org

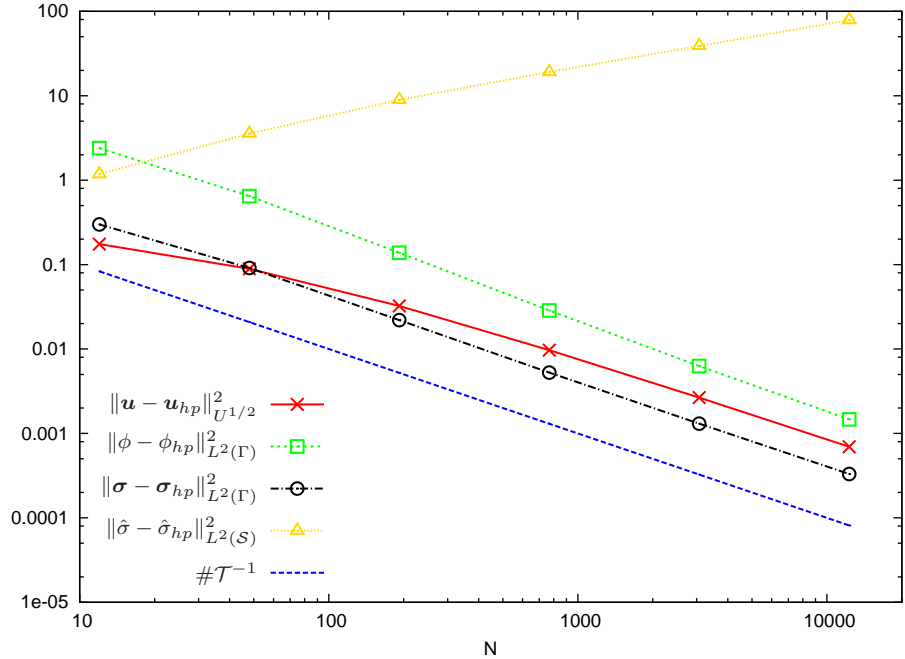

Figure 1: Experiment I. Errors for smooth solution on a closed surface.

**Experiment II.** The surface $\Gamma$ and the coarsest mesh $\mathcal{T}_0$ are chosen as in Experiment I. The right-hand side is chosen as $f(x,y,z) = x$. Sine $\Omega$ is centered at the origin, the mean value of $f$ on $\Gamma$ vanishes. In this case there holds $\phi, f \in H^1(\Gamma)$ and $\boldsymbol{\sigma} \in \mathbf{H}^1(\Gamma)$. The resulting convergence order $O(h)$ for uniform mesh-refinement is confirmed by Figure 2.

**Experiment III.** We choose $\Gamma$ to be the square $(-1,1)^2$. The coarsest mesh consists of 4 triangles that appear when $\Gamma$ is divided by its 2 diagonals. The exact solution $\phi$ is prescribed as a piecewise affine, globally continuous function, with value 1 at the node in the center of $\Gamma$, and values 0 at the nodes on the boundary of $\Gamma$. The outcome of our method with uniform mesh-refinement is shown in Figure 3 and indicates convergence order $O(h)$. The comments from Experiment I regarding the regularity of $\hat{\sigma}$ apply. For illustration we show some sample solutions $\phi_{hp}$ in Figure 4. They perfectly approximate the selected exact solution $\phi$ by piecewise constants.

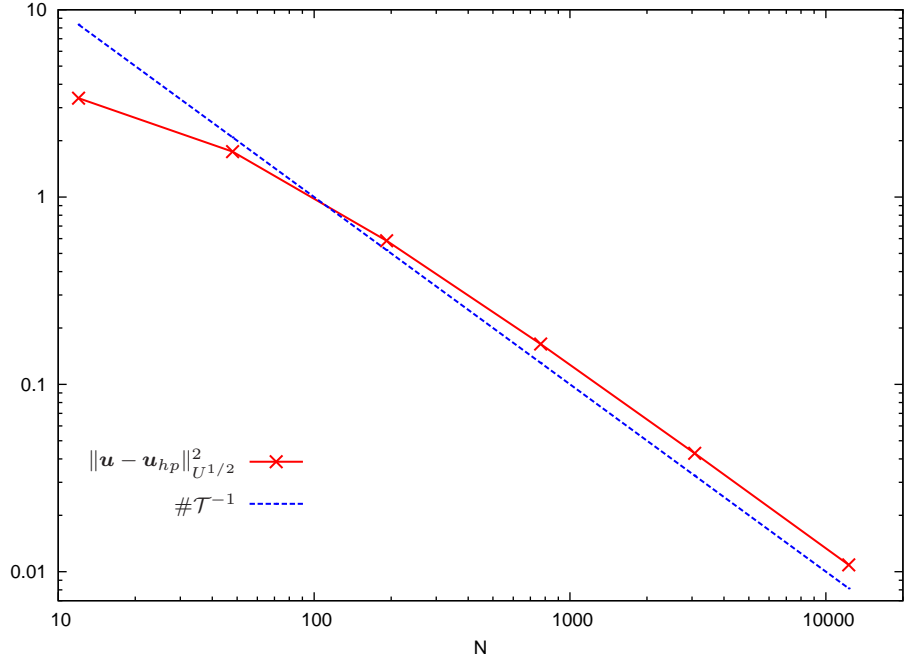

Figure 2: Experiment II. Error for singular solution on a closed surface.

**Experiment IV.** The surface $\Gamma$ and the coarsest mesh $\mathcal{T}_0$ are chosen as in Experiment III. The right-hand side is given by $f(x,y,z) = 1$. This is the case where $\phi$ is expected to have strong edge singularities so that $\phi \in H^t(\Gamma)$ for any $t < 1$, but $t \neq 1$ in general, cf. [38]. Figure 5 shows the error in energy norm (squared) for uniform mesh-refinement. It confirms the order $O(h)$, which is the limit case (it is just excluded from the theory). Also for this experiment we show some sample solutions $\phi_{hp}$ (Figure 6). They resemble the typical shape of $\phi$ with square root singularity $\mathrm{dist}_\Gamma(\cdot, \partial\Gamma)^{1/2}$.

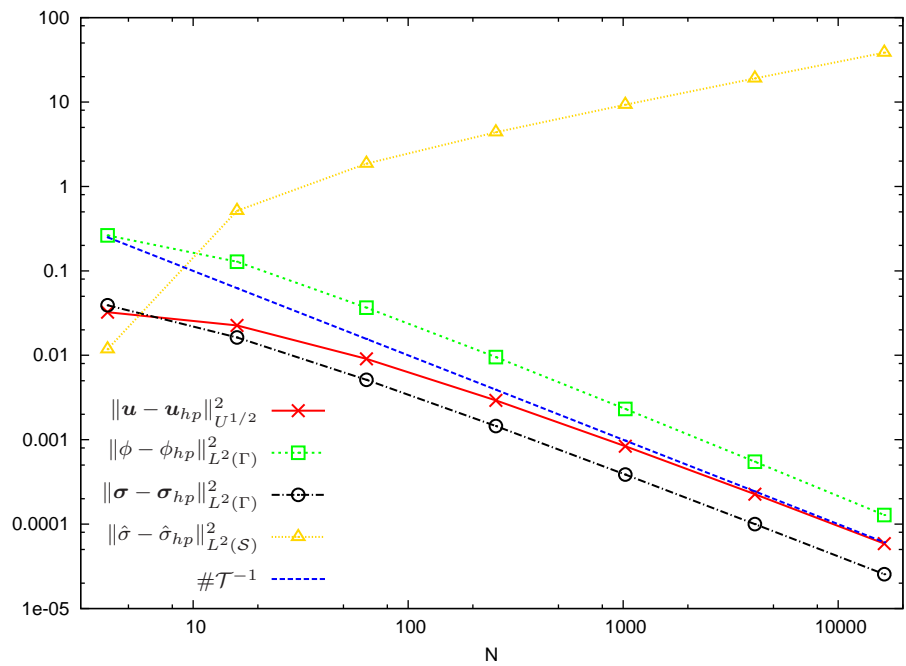

Figure 3: Experiment III. Errors for smooth solution on an open surface.

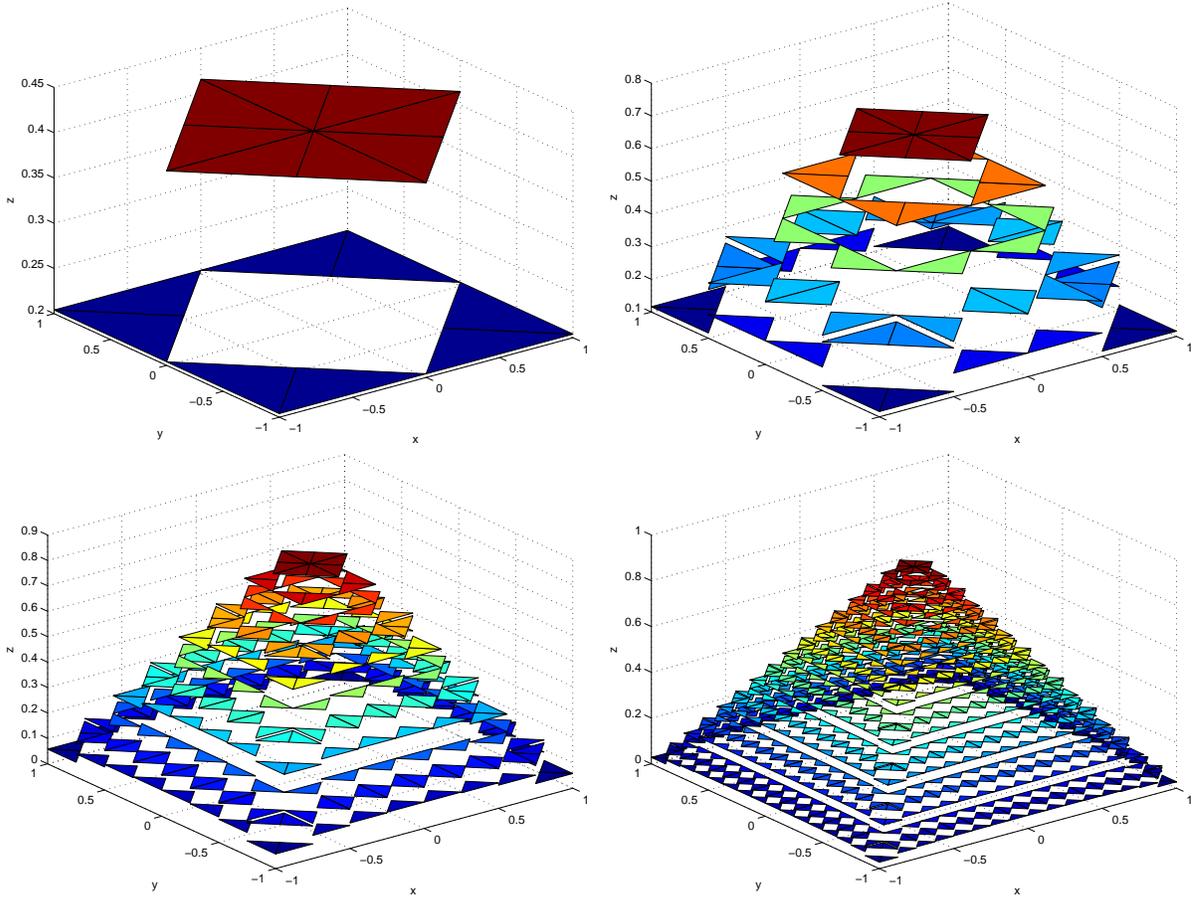

Figure 4: Solutions $\phi_{hp}$ of Experiment III on the meshes $\mathcal{T}_i$ for $i = 1, \ldots 4$.

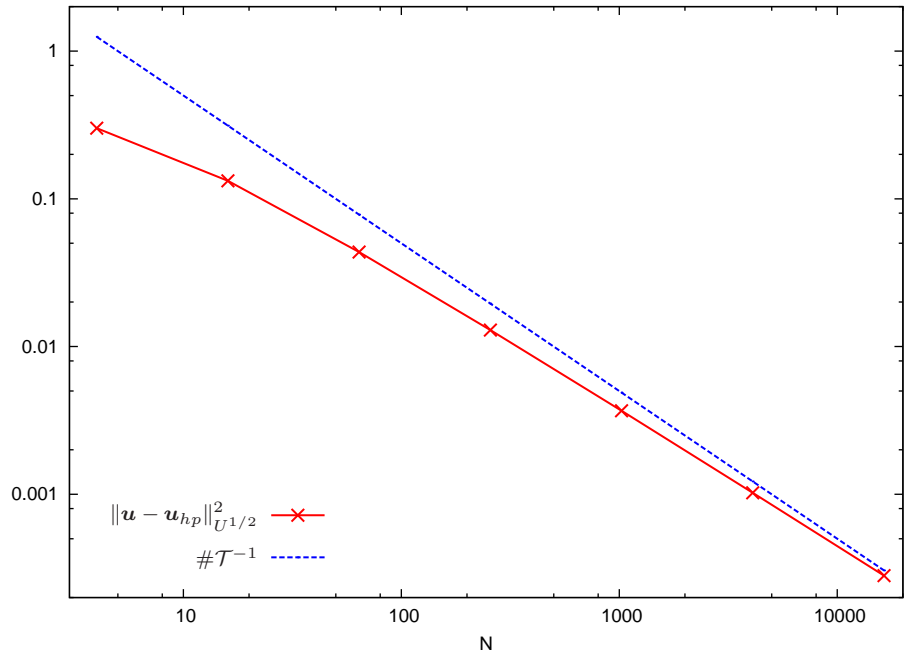

Figure 5: Experiment IV. Error for singular solution on an open surface.

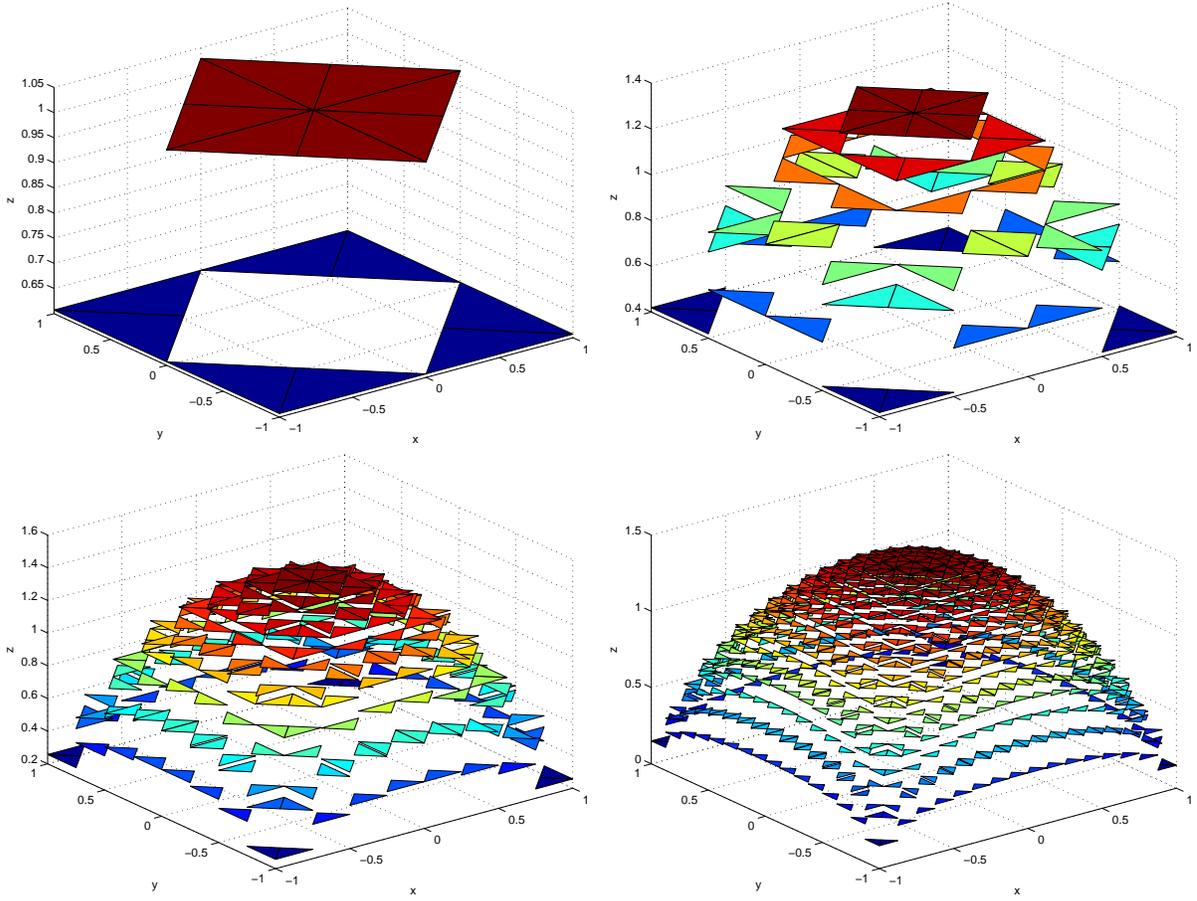

Figure 6: Solutions $\phi_{hp}$ of Experiment IV on the meshes $\mathcal{T}_i$ for $i = 1, \ldots 4$.